\documentclass[a4paper,10pt]{article}
\usepackage{graphicx}
\usepackage{amsmath}
\usepackage{amssymb}
\usepackage{subfig}
\usepackage{booktabs}
\usepackage{multirow}
\usepackage{tikz}
\usepackage{graphicx}
\usetikzlibrary{shapes,arrows}
\usepackage{sidecap}

\newtheorem{theorem}{Theorem}[section]
\newtheorem{lemma}[theorem]{Lemma}
\newtheorem{definition}[theorem]{Definition}
\newcommand{\vertiii}[1]{{\left\vert\kern-0.25ex\left\vert\kern-0.25ex\left\vert #1 
    \right\vert\kern-0.25ex\right\vert\kern-0.25ex\right\vert}}

\title{High order finite difference methods for the wave equation with non--conforming grid interfaces}

\author{Siyang Wang\footnotemark[1]\ \footnotemark[2]
\and Kristoffer Virta\footnotemark[2]
\and Gunilla Kreiss\footnotemark[2]}

\date{}
\begin{document}
\maketitle

\renewcommand{\thefootnote}{\fnsymbol{footnote}} 
\footnotetext[1]{Corresponding author, email: {siyang.wang@it.uu.se}}
 \footnotetext[2]{Department of Information Technology, Uppsala University, Uppsala, SE-751 05, Sweden}
\renewcommand{\thefootnote}{\arabic{footnote}}

\begin{abstract}
We use high order finite difference methods to solve the wave equation in the second order form. The spatial discretization is performed by finite difference operators satisfying a summation--by--parts property. The focus of this work is on the numerical treatment of non--conforming grid interfaces. The interface conditions are imposed weakly by the simultaneous approximation term technique in combination with interface operators, which move the discrete solutions between the grids on the interface. In particular, we consider interpolation operators and projection operators. A norm--compatibility condition, which leads to stability for first order hyperbolic systems, does not suffice for second order wave equations. An extra constraint on the interface operators must be satisfied to derive an energy estimate for stability. We carry out eigenvalue analyses to investigate the additional constraint and how it is related to stability, and find that the projection operators have better stability properties than the interpolation operators. In addition, a truncation error analysis is performed to study the convergence property of the numerical schemes. In the numerical experiments, the stability and accuracy properties of the numerical schemes are further explored, and the practical usefulness of non--conforming grid interfaces is presented and discussed in two efficiency studies.
\end{abstract}



\pagestyle{myheadings}
\thispagestyle{plain}

\section{Introduction}
For wave propagation problems, the computational domain is often large compared with the wavelength, and waves travel for a long time. It has been shown that high order accurate discretizations in time and space are more efficient to solve these problems on smooth domains \cite{Gustafsson2008,Kreiss1972}. Although it is straightforward to derive high order finite difference schemes in the interior of the computational domain, it is challenging to derive the boundary closures in a stable and accurate way. For long time simulations, it is also desirable that the discretization is strictly stable \cite[pp.~129]{Gustafsson2013}. A successful candidate of high order finite differences is the summation--by--parts simultaneously approximation term (SBP--SAT) method \cite{Fernandez2014,Svard2014}. An SBP operator \cite{Kreiss1974} approximates a spatial derivative, and mimics integration--by--parts via its associated norm. The SAT method \cite{Carpenter1994} is used to impose boundary conditions and grid interface conditions weakly.  

Traditionally, the wave equation is written as a first order hyperbolic system, and is then solved by the well--developed methods for such systems. However, there are various drawbacks in doing so \cite{Kreiss2002}. Therefore, it is desirable to solve directly the wave equation in the second order form. In \cite{Virta2014}, an SBP--SAT method for the wave equation in the second order form is developed and the numerical treatment of conforming grid interfaces is discussed. Stability of the numerical scheme is proved by the energy method, and the convergence property is investigated in the numerical experiments.

For a wave that travels in an inhomogeneous medium, the wave speed varies in space. Since the wavelength is proportional to the wave speed, a reduction in the wave speed confined to a subset of the physical domain yields a wave with a shorter wavelength localized in that subset. In \cite{Kreiss1972}, the accuracy of a numerical solution to a Cauchy problem is stated in terms of the number of grid points per wavelength. For computational efficiency it is important that a fine mesh is used in the subset that constitutes the slower media, and a coarse mesh elsewhere. 

To achieve this, one approach is to partition the computational domain into blocks, where the mesh sizes are constant in each block but differ in different blocks. In more than one space dimension, the partition results in non--conforming grid interfaces with hanging nodes. Suitable interface conditions are then imposed to couple adjacent mesh blocks and yield a well--posed problem. Many techniques for the numerical treatment of interface conditions have been proposed. In \cite{Petersson2010}, an energy conserving discretization of the elastic wave equation in the second order formulation is presented. The finite difference operators satisfy the principle of SBP, and the grid interface with a 1:2 refinement ratio is handled by ghost points.  Stability is proved by the energy method, but the convergence rate of the numerical scheme is limited to two. 

In \cite{Mattsson2010} the norm--compatible interpolation operators are constructed to handle non--conforming grid interfaces, and the Euler equations are used as the model problem. The norm--compatibility condition leads to an energy estimate for first order hyperbolic systems and the Schr\"{o}dinger equation \cite{Nissen2012}. The interpolation error is of the same magnitude as the error due to the derivative approximations by the SBP operators. In the numerical experiments, it is shown that the convergence rate of the numerical scheme applied to the Euler equations is not lowered by using the interpolation operators, compared with the case with only conforming grid interfaces. To use the interpolation operators presented in \cite{Mattsson2010}, the mesh refinement ratio is fixed to $1:2$ and the mesh blocks must be conforming.  

Recently, a general purpose methodology for coupling mesh blocks with non--conforming interfaces was developed in \cite{Kozdon2015}. This technique uses projection operators to move a discrete finite difference solution to piecewise polynomial functions in a subspace of a Hilbert space where the coupling is done. The wave equation in the first order system formulation is used as the model problem. Stability is proved by the energy method, and the numerical experiments demonstrate that the convergence rate is the same as if conforming grid interfaces were used. The projection operators allow for a very flexible configuration of meshing in the sense that the interfaces as well as the mesh blocks do not need to be conforming. Similarly to the interpolation operators, the projection operators satisfy norm--compatibility conditions that are essential for the stability proof to hold. 

In this paper, we focus on the numerical treatment of non--conforming grid interfaces for the wave equation in the second order form in the framework of the SBP--SAT methodology. In particular, the stability and accuracy properties are investigated. We have found that in contrast to first order hyperbolic systems for which the norm--compatibility condition leads to stability, an extra condition on the interface operator is needed to derive an energy estimate for the second order wave equation. This condition is satisfied for the second and fourth order accurate interpolation operators constructed in \cite{Mattsson2012} and the projection operators in \cite{Kozdon2015}. We prove stability by the energy method for those cases.  For higher order accurate schemes the extra condition is not satisfied, and we cannot prove stability. With an eigenvalue analysis, we have found that the violation of the stability condition is very weak for the projection operators, and in all the numerical experiments we have conducted no unphysical growth is observed for the schemes with the projection operators. In certain cases the SBP--SAT schemes with the sixth and eighth order accurate interpolation operators are not stable. 

Local mesh refinement reduces the number of grid points significantly in computations. To achieve full efficiency, the numerical scheme must also be accurate enough. It is desirable that the convergence rate is not depressed by using non--conforming grid interfaces. Even though this is in most cases true for first order hyperbolic systems, the situation for second order equations is less favourable. By a truncation error analysis, we show that the truncation error near the edge of the non--conforming grid interfaces is two orders larger than that with conforming grid interfaces. The large truncation error is only localized at a few grid points in a two dimensional space, and its effect to the convergence rate may be weakened. In fact, the numerical experiments show that the convergence rate with a non--conforming grid interface is only one order lower than the corresponding case with a conforming grid interface. In addition, an efficiency study is carried out by a comparison of the numerical schemes with interpolation operators and projection operators. We have found that in certain cases it is beneficial to use non--conforming grid interfaces, albeit the accuracy reduction.

The structure of this paper is as follows. In \S\ref{sec-Preliminaries}, the SBP--SAT methodology is introduced. We then discuss stability and accuracy properties of the numerical coupling based on interpolation operators in \S\ref{sec-MC}, and projection operators in \S\ref{sec-KW}. Numerical experiments are carried out in \S\ref{sec-NE} including the eigenvalue analyses for stability, convergence verifications for accuracy and studies on computational efficiency. We  conclude and mention future work in \S\ref{sec-Conclusion}. 

\section{Preliminaries}\label{sec-Preliminaries}
We begin with the preliminaries that will be used in the discussion of the SBP--SAT method. Let $w_1(x)$ and $w_2(x)$ be two real--valued functions in L$_2[0,1]$. The inner product is defined by $(w_1,w_2)=\int_0^1w_1w_2dx$. The corresponding norm is $\|w_1\|^2=(w_1,w_1)$. The computational domain $[0,1]$ is discretized by $N+1$ equidistant grid points
\begin{equation*}
x_i=ih,\quad i=0,1,\cdots,N, \text{ where } h=\frac{1}{N}.
\end{equation*}
With any fixed $N$, a grid function can be represented by a vector and an operator can be represented by a matrix. Throughout this paper, we use an operator and a matrix interchangeably when there is no ambiguity.  

\subsection{The SBP operators}
We need the following definitions:
\begin{definition}
A difference operator $D_1=H^{-1}Q$ approximating first derivative $\partial/\partial x$ is a narrow diagonal first derivative SBP operator, if $H$ is diagonal and positive definite, $Q+Q^T = B = \text{diag} (-1, 0,... , 0, 1)$ and the interior stencil  width of $D_1$ is minimal.
\end{definition}
\begin{definition}\label{defD2}
A difference operator $D_2^{(b)}=H^{-1}(-M^{(b)}+B^{(b)}S)$ approximating second derivative $\partial /\partial x (b(x)\partial /\partial x)$ with $b(x)>0$ is a narrow diagonal variable coefficient second derivative SBP operator, if $H$ is diagonal and positive definite, $M^{(b)}$ is symmetric positive semidefinite, $B^{(b)} = \text{diag} (-b(x_0), 0,... , 0, b(x_N))$, the first and last row of $S$ is a one sided approximation of $\partial/\partial x$ at the boundary and the interior stencil width of $D_2^{(b)}$ is minimal.
\end{definition}

The diagonal positive definite matrix $H$ defines the SBP norm, and it has the interior weight $h$ and special boundary weights. To solve the wave equation on a curvilinear grid, even if the original equation has only constant coefficients, the transformed equation has second derivative terms with variable coefficients, and mixed derivative terms. Therefore, it is important that $D_1$ and $D_2^{(b)}$ are based on the same norm $H$. In addition, $D_1$ and $D_2^{(b)}$ must be compatible for the energy method to be applicable for proving stability. Compatibility \cite{Mattsson2008b} means that   $M^{(b)}$ can be written as $M^{(b)}=D_1^THB^{(b)}D_1+R^{(b)}$ where $R^{(b)}$ is symmetric positive semidefinite. In this case, $D_2^{(b)}$ is essentially equal to applying $D_1$ twice plus a small dissipative term.

The accuracy of the SBP operators are often termed as $2p$, meaning that the approximation error of $D_1$ and $D_2^{(b)}$ is $\mathcal{O}(h^{2p})$ in the interior, while near the boundary the approximation error increases to $\mathcal{O}(h^{p})$. For $S\approx \frac{d}{dx}$ at the boundary, the approximation error is $\mathcal{O}(h^{p+1})$. $2p^{th}$ order accurate SBP operators $D_1$ are constructed in \cite{Strand1994} for $p=1,2,3,4$, and in \cite{Mattsson2013} for $p=5$. In \cite{Mattsson2012}, $2p^{th}$ order accurate $D_2^{(b)}$ are constructed for $p=1,2,3$, and they are compatible with $D_1$ constructed in \cite{Strand1994}. The construction of $D_2^{(b)}$ requires to solve a large system of nonlinear equations, which makes it very involved when the accuracy order is high. In certain cases, we only need an SBP operator $D_2$ that approximates second derivative $\partial^2 /\partial x^2$, and these operators are constructed in \cite{Mattsson2004} for $p=1,2,3,4$ and in \cite{Mattsson2013} for $p=5$.

The following lemma, which is often referred to as the \emph{borrowing trick} \cite{Mattsson2008}, describes an important property of the SBP operator $D_2$ constructed in \cite{Mattsson2013,Mattsson2004}:
\begin{lemma}\label{borrow}
The matrix $M$ in $D_2$ constructed in \cite{Mattsson2004} can be written as
\begin{equation*}
M=h\alpha_{2p}(BS)^TBS+\tilde{M},
\end{equation*}
where $\tilde{M}$ is symmetric positive semidefinite and $\alpha_{2p}$ is a constant independent of $h$. The values of $\alpha_{2p}$ are listed in Table \ref{alpha} for different accuracy orders. 
\end{lemma} 
\begin{table}
\centering
\begin{tabular}{c c  c  c c c}
\toprule
$2p$ & 2 & 4 & 6 & 8 & 10\\ \midrule
$\alpha_{2p}$ &  0.4 & 0.2508560249 & 0.1878715026 & 0.0015782259 & 0.0351202265\\ \bottomrule
\end{tabular}
\caption{$\alpha_{2p}$ values}
\label{alpha}
\end{table}

A similar lemma for $D_2^{(b)}$ constructed in \cite{Mattsson2012} is presented in \cite{Virta2014}.

There are other types of SBP operators as well. It is possible to increase the accuracy by using a block norm $H$, meaning that $H$ is diagonal in the interior and has a symmetric--block--structured boundary closure. In this case, the approximation error in the interior and near the boundary are $\mathcal{O}(h^{2p})$ and $\mathcal{O}(h^{2p-1})$, respectively. The drawback of the block norm SBP operators is that the energy method is not applicable to prove stability for the wave equation with variable coefficients and/or solved on a curvilinear grid, and unphysical growth is indeed observed in numerical experiments. This has limited the practical usage of the block norm SBP operators. In \cite{Mattsson2013}, artificial dissipation is added to stablize the scheme.

\section{Non--conforming grid interfaces handled by interpolation operators}\label{sec-MC}
\begin{figure}
\begin{center}
\subfloat[Two mesh blocks]{\includegraphics[width=4cm]{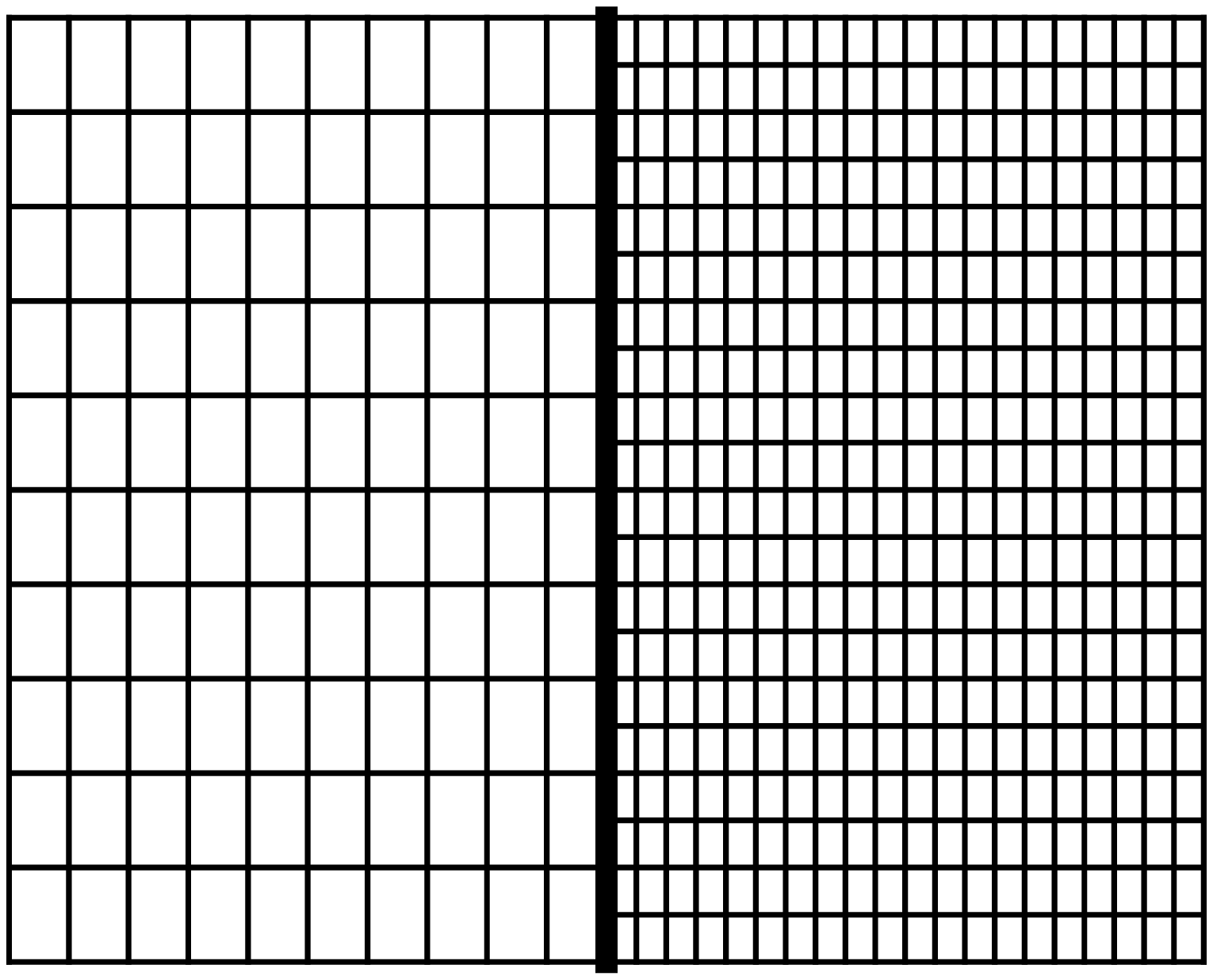}\label{Mesh_2_MC}}
\subfloat[Four mesh blocks]{\includegraphics[width=4cm]{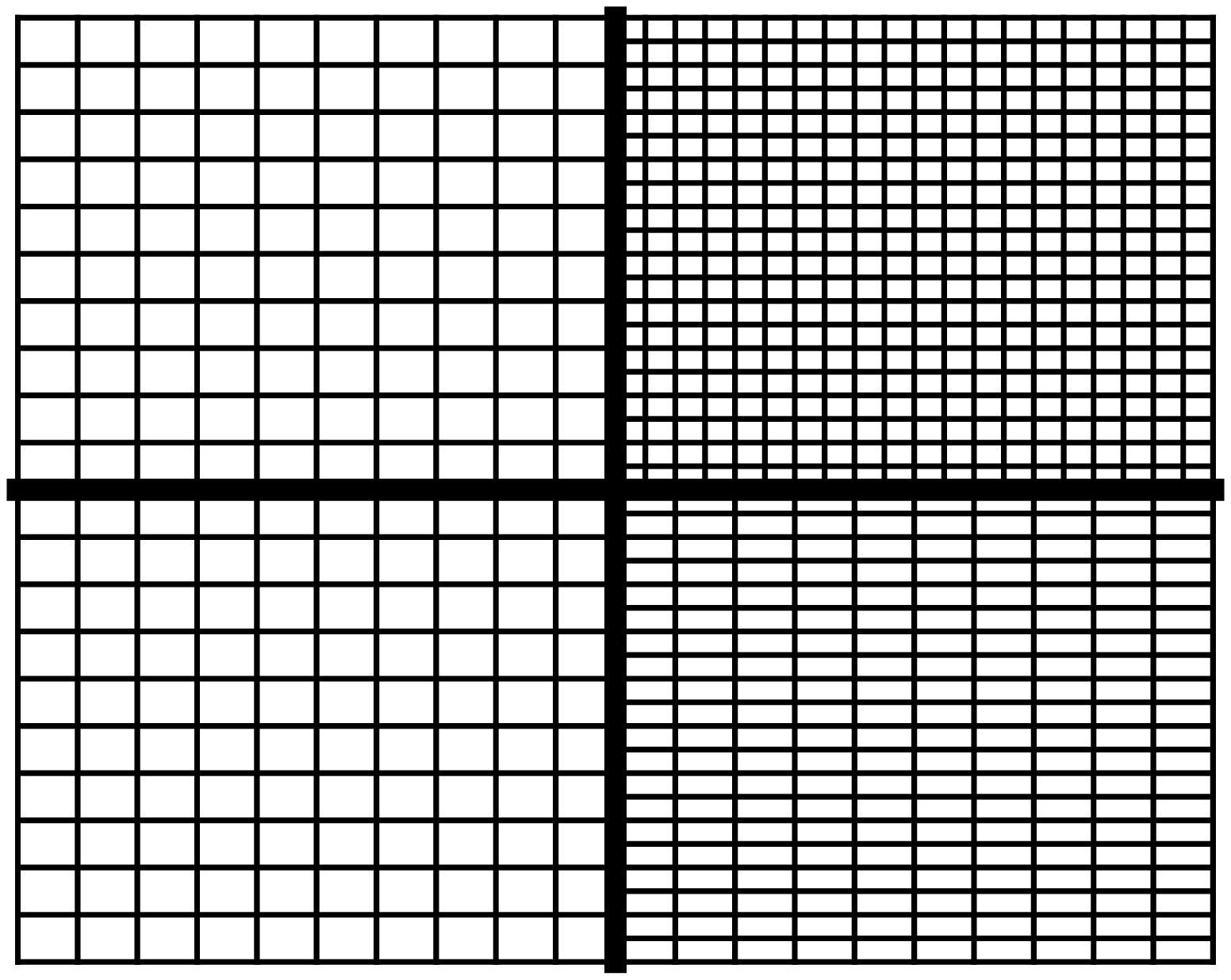}\label{Mesh_4_MC}}
\caption{Non--conforming grid interfaces}
\label{Mesh_MC}
\end{center}
\end{figure}

The numerical coupling of conforming grid interfaces by the SAT method for the wave equation is discussed in \cite{Virta2014}. Our aim in this section is to generalize the scheme to also couple non--conforming grid interfaces by using the interpolation operators constructed in \cite{Mattsson2010}. To begin with, we consider the two--block--structured mesh $\Omega$ shown in Figure \ref{Mesh_2_MC}: a coarse mesh $\Omega_L$ in the left block with $n_{xL}\times n_{yL}$ grid points and a fine mesh $\Omega_R$ in the right block with $n_{xR}\times n_{yR}$ grid points. The equality $n_{yR}=2n_{yL}-1$ yields a 1:2 mesh refinement ratio across the interface. 

In \cite{Mattsson2010}, the interpolation operators are denoted by $I_{F2C}$ and $I_{C2F}$, where the subscripts $F2C$ and $C2F$ refer to \emph{fine to coarse} and \emph{coarse to fine}, respectively. Though only interpolation operators for the grid interface with mesh refinement ratio 1:2 are reported, the construction of the interpolation operators with some other refinement ratios (1:4, 1:8, $\cdots$) is possible by the same technique. To handle a multi--block--structured mesh shown in Figure \ref{Mesh_4_MC}, even though the interface between the left--up block and the left--down block is conforming, it is necessary to treat it as an interface to make the energy method applicable to prove stability. In other words, the mesh blocks must be conforming.

It is important that the interpolation operators preserve the SBP property. To this end, the following norm--compatibility condition is essential:
\begin{equation}\label{equality}
H_{yR}I_{C2F}=(H_{yL}I_{F2C})^T,
\end{equation}
where $H_{yL}$ and $H_{yR}$ are SBP norms in the left block and right block, both in the $y$ direction. 

The interpolation operators do not interpolate exactly, instead they mimic the accuracy properties of the diagonal norm SBP operators. The interpolation error is $\mathcal{O}(h^{2p})$ in the interior of the interface and $\mathcal{O}(h^{p})$ near the edge of the interface. We call them $2p^{th}$ order accurate interpolation operators, and when used together with the $2p^{th}$ order accurate SBP operators the scheme is also termed as $2p^{th}$ order accurate, though the truncation error of the semidiscretized equation may not be $\mathcal{O}(h^{2p})$ or $\mathcal{O}(h^{p})$. Here, $h$ is used to denote the magnitude of the mesh sizes for the sake of a simplified notation, though at most four different mesh sizes could be present in the mesh $\Omega$.

With the above accuracy requirement and the norm--compatibility condition (\ref{equality}), $2p^{th}$ order accurate interpolation operators are constructed for $p=1,2,3$ and 4 in \cite{Mattsson2010}. Though the accuracy is reduced near the edge of the interface, the number of grid points with the large interpolation error $\mathcal{O}(h^p)$ is independent of $h$. 

As will be seen later, when the interpolation operators are used to solve the wave equation an extra condition is posed on the interpolation operators in order to apply the energy method to prove stability:
\begin{equation}\label{inequality}
\Xi_L:=H_{yL}(I_{yL}-I_{F2C}I_{C2F})\geq 0,\quad \Xi_R:=H_{yR}(I_{yR}-I_{C2F}I_{F2C})\geq 0,
\end{equation}
where $\Xi\geq 0$ means that the matrix $\Xi$ is symmetric positive semidefinite, and $I_{yL}$, $I_{yR}$ are identity matrices. It is straightforward to show that $\Xi_L$ and $\Xi_R$ are symmetric, but the positive semidefiniteness is not a built--in constraint in the construction process of the interpolation operators in \cite{Mattsson2010}. In \S\ref{sec-NE}, an eigenvalue analysis is performed to show that (\ref{inequality}) is satisfied for the second and fourth order accurate interpolation operators. Negative eigenvalues are present with the sixth and eighth order accurate interpolation operators. However, the scheme with sixth other accurate interpolation operators seems stable in some of the numerical experiments conducted in \S\ref{sec-NE}, indicating that  (\ref{inequality}) is a sufficient but not necessary condition for stability.

\subsection{The wave equation with a non--conforming grid interface}
The wave equation in the second order form in two space dimensions is
\begin{equation}\label{Wave2d}
U_{tt}=U_{xx}+U_{yy},
\end{equation}
where $-\infty<x<\infty$, $0\leq y\leq 1$ and $0\leq t\leq t_f$. We assume that the initial conditions and boundary conditions are compatible smooth functions with compact support. As a consequence, the true solution $U$ is also smooth. To solve the equation on the mesh $\Omega$ shown in Figure \ref{Mesh_2_MC}, continuity of the solution and continuity of first normal derivative across the grid interface are required. 

In the numerical coupling scheme, we frequently pick up solutions along the interface by using the matrices defined in the first column of Table \ref{TableM}. In each of those matrices, all elements are zero except one element that is equal to one. The sizes along with the positions of the nonzero element are listed in the second and third column of Table \ref{TableM}. Note that $E_{LR}=E_{RL}^T$.
\begin{SCtable}[.5]
\caption{Matrices that are used to pick up solutions on the interface.}
\begin{tabular}{c c c}
\toprule
Matrix & Size & Nonzero\\ \midrule
$E_{0L}$ & $n_{xL}\times n_{xL}$ & $(n_{xL}, n_{xL})$\\ \midrule
$E_{0R}$ & $n_{xR}\times n_{xR}$ & $(1,1)$\\ \midrule
$E_{LR}$ & $n_{xL}\times n_{xR}$ & $(n_{xL}, 1)$\\ \midrule
$E_{RL}$ & $n_{xR}\times n_{xL}$ & $(1, n_{xR})$\\
\bottomrule
\end{tabular}
\label{TableM}
\end{SCtable}

Bold letters are used to denote the operators in two space dimensions, which are obtained from the corresponding one dimensional operators through the Kronecker product $\boldsymbol{A_x}=A_x\otimes I_y$ and $\boldsymbol{A_y}=I_x\otimes A_y$, where $I_x$ and $I_y$ are identity matrices. 

Next, Equation (\ref{Wave2d}) is discretized by the SBP operators on $\Omega$ and the grid interface conditions are imposed weakly by the SAT method and interpolation operators. The semidiscretized equation corresponding to (\ref{Wave2d}) reads:
\begin{subequations}\label{Wave2dsemi}
\begin{align}
&\boldsymbol{u_{tt}=D_{2L}u+SAT_{u1}+SAT_{u2}+SAT_{\partial u}}\label{Wave2dsemiu}, \\
&\boldsymbol{v_{tt}=D_{2R}v+SAT_{v1}+SAT_{v2}+SAT_{\partial v}}\label{Wave2dsemiv},
\end{align}
\end{subequations}
where
\begin{equation*}
\begin{split}
&\boldsymbol{SAT_{u1}}=\frac{1}{2}\boldsymbol{H_{xL}^{-1}S_{xL}^T(E_{0L}u-}(E_{LR}\otimes I_{F2C}) \boldsymbol{v}), \\
&\boldsymbol{SAT_{u2}}=-\tau\boldsymbol{H_{xL}^{-1}(E_{0L}u-}(E_{LR}\otimes I_{F2C}) \boldsymbol{v}), \\
&\boldsymbol{SAT_{\partial u}}=-\frac{1}{2}\boldsymbol{H_{xL}^{-1}(E_{0L}S_{xL}u-} (E_{LR}\otimes I_{F2C}) \boldsymbol{S_{xR}v}),
\end{split}
\end{equation*}
and 
\begin{equation*}
\begin{split}
&\boldsymbol{SAT_{v1}}=-\frac{1}{2}\boldsymbol{H_{xR}^{-1}S_{xR}^T(E_{0R}v-}(E_{RL}\otimes I_{C2F}) \boldsymbol{u}), \\
&\boldsymbol{SAT_{v2}}=-\tau\boldsymbol{H_{xR}^{-1}(E_{0R}v-}(E_{RL}\otimes I_{C2F}) \boldsymbol{u}), \\
&\boldsymbol{SAT_{\partial v}}=\frac{1}{2}\boldsymbol{H_{xR}^{-1}(E_{0R}S_{xR}v-} (E_{RL}\otimes I_{C2F}) \boldsymbol{S_{xL}u}).
\end{split}
\end{equation*}
Here, $\boldsymbol{u}$ and $\boldsymbol{v}$ are grid functions in $\Omega_L$ and $\Omega_R$, respectively. $\boldsymbol{D_{2L}}$ and $\boldsymbol{D_{2R}}$ are SBP operators approximating second derivatives. In (\ref{Wave2dsemi}) the penalty terms for the boundary conditions  are omitted as the focus here is the numerical treatment of the grid interface. Both $\boldsymbol{SAT_{u1}}$ and $\boldsymbol{SAT_{u2}}$ impose weakly the continuity of the solution across the grid interface. The term $\boldsymbol{S_{xL}^T}$ in $\boldsymbol{SAT_{u1}}$ makes the semidiscretization symmetric with respect to the SBP norms. The penalty parameter $\tau$ in $\boldsymbol{SAT_{u2}}$ controls the strength of the weak enforcement, and its value is determined by the energy method. The penalty term $\boldsymbol{SAT_{\partial u}}$ imposes weakly the continuity of the first normal derivative across the grid interface. The penalty terms in Equation (\ref{Wave2dsemiv}) are constructed in a similar way.

In the following, the stability of (\ref{Wave2dsemi}) is proved by the energy method.

\subsubsection{Stability analysis by the energy method}
Multiplying Equation (\ref{Wave2dsemiu}) by $\boldsymbol{u_t^TH_L}$ and Equation (\ref{Wave2dsemiv}) by $\boldsymbol{v_t^TH_R}$, and using the equality $E_{LR}=E_{RL}^T$ and relation (\ref{equality}), 
\begin{equation}\label{E2d1}
\begin{split}
&\frac{d}{dt}(\boldsymbol{u_t^TH_Lu_t+v_t^TH_Rv_t}) \\
=&\frac{d}{dt}(\boldsymbol{-u^T}(M_{xL}\otimes H_{yL})\boldsymbol{u} + \boldsymbol{u^T}(E_{0L}S_{xL}\otimes H_{yL})\boldsymbol{u} - \tau \boldsymbol{u^T}(E_{0L}\otimes H_{yL}) \boldsymbol{u} \\
& \boldsymbol{-v^T}(M_{xR}\otimes H_{yR})\boldsymbol{v} - \boldsymbol{v^T}(E_{0R}S_{xR}\otimes H_{yR})\boldsymbol{v} - \tau \boldsymbol{v^T}(E_{0R}\otimes H_{yR})\boldsymbol{v} \\
&-\boldsymbol{u^T}(S_{xL}^T E_{LR}\otimes H_{yL} I_{F2C})\boldsymbol{v} + \boldsymbol{v^T}(S_{xR}^TE_{RL}\otimes H_{yR}I_{C2F})\boldsymbol{u} \\
&+ 2\tau \boldsymbol{u^T}(E_{LR}\otimes H_{yL}I_{F2C})\boldsymbol{v}).
\end{split}
\end{equation}
Next, the following equality is obtained by moving all terms on the right hand side of (\ref{E2d1}) to the left
\begin{equation}
\frac{d}{dt}\boldsymbol{E_H^W}=0,
\end{equation}
where 
\begin{equation}\label{EHW}
\begin{split}
\boldsymbol{E_H^W}=&\boldsymbol{u_t^TH_Lu_t+v_t^TH_Rv_t} \\
&+\boldsymbol{u^T}(M_{xL}\otimes H_{yL})\boldsymbol{u} - \boldsymbol{u^T}(E_{0L}S_{xL}\otimes H_{yL})\boldsymbol{u} + \tau \boldsymbol{u^T}(E_{0L}\otimes H_{yL}) \boldsymbol{u} \\
& +\boldsymbol{v^T}(M_{xR}\otimes H_{yR})\boldsymbol{v} + \boldsymbol{v^T}(E_{0R}S_{xR}\otimes H_{yR})\boldsymbol{v} + \tau \boldsymbol{v^T}(E_{0R}\otimes H_{yR})\boldsymbol{v} \\
&+\boldsymbol{u^T}(S_{xL}^T E_{LR}\otimes H_{yL} I_{F2C})\boldsymbol{v} - \boldsymbol{v^T}(S_{xR}^TE_{RL}\otimes H_{yR}I_{C2F})\boldsymbol{u} \\
&-2\tau \boldsymbol{u^T}(E_{LR}\otimes H_{yL}I_{F2C})\boldsymbol{v}.
\end{split}
\end{equation}
Then, an energy estimate exists if $\boldsymbol{E_H^W}\geq 0$ for any $\boldsymbol{u}$ and $\boldsymbol{v}$, and in this case $\boldsymbol{E_H^W}$ is a discrete energy of the semidiscretized equation (\ref{Wave2dsemi}). To achieve this, we need relation (\ref{inequality}) to obtain  
\begin{equation*}
\begin{split}
\boldsymbol{u^T}(E_{0L}\otimes H_{yL}) \boldsymbol{u} &\geq\boldsymbol{u^T}(E_{0L}\otimes H_{yL}I_{F2C}I_{C2F}) \boldsymbol{u} \\
&=\boldsymbol{u^T}(E_{0L}\otimes I_{C2F}^TH_{yR}I_{C2F}) \boldsymbol{u} \\
&=((E_{0L}\otimes I_{C2F})\boldsymbol{u})^T (I_{xL}\otimes H_{yR}) ((E_{0L}\otimes I_{C2F})\boldsymbol{u}).
\end{split}
\end{equation*}
It then follows for $\boldsymbol{u}$:
\begin{equation}\label{trick2u}
\boldsymbol{u^T}(E_{0L}\otimes H_{yL}) \boldsymbol{u} \geq \frac{1}{2}\boldsymbol{u^T}(E_{0L}\otimes H_{yL}) \boldsymbol{u} + \frac{1}{2}((E_{0L}\otimes I_{C2F})\boldsymbol{u})^T (I_{xL}\otimes H_{yR}) ((E_{0L}\otimes I_{C2F})\boldsymbol{u}),
\end{equation}
and for $\boldsymbol{v}$:
\begin{equation}\label{trick2v}
\boldsymbol{v^T}(E_{0R}\otimes H_{yR}) \boldsymbol{v} \geq \frac{1}{2}\boldsymbol{v^T}(E_{0R}\otimes H_{yR}) \boldsymbol{v} + \frac{1}{2}((E_{0R}\otimes I_{F2C})\boldsymbol{v})^T (I_{xR}\otimes H_{yL}) ((E_{0R}\otimes I_{F2C})\boldsymbol{v}).
\end{equation}
In addition, Lemma \ref{borrow} gives:
\begin{equation}\label{borrowM}
\begin{split} 
&M_{xL}=\tilde{M}_{xL}+h_{xL}\alpha_{2p}(E_{0L}S_{xL})^T (E_{0L}S_{xL}), \\
&M_{xR}=\tilde{M}_{xR}+h_{xR}\alpha_{2p}(E_{0R}S_{xR})^T (E_{0R}S_{xR}),
\end{split}
\end{equation}
where $\tilde{M}_{xL}$ and $\tilde{M}_{xR}$ are symmetric positive semidefinite matrices. 

Next, we plug in (\ref{borrow}), (\ref{trick2u}) and (\ref{trick2v}) to ({\ref{EHW}), and obtain
\begin{equation*}
\boldsymbol{E_H^W=Q_1+Q_2+Q_3},
\end{equation*}
where 
\begin{equation*}
\begin{split}
\boldsymbol{Q_1} =& \boldsymbol{u_t^TH_Lu_t+v_t^TH_Rv_t}+\boldsymbol{u^T}(\tilde{M}_{xL}\otimes H_{yL})\boldsymbol{u}+ \boldsymbol{v^T}(\tilde{M}_{xR}\otimes H_{yR})\boldsymbol{v}, \\
\boldsymbol{Q_2} =& h_{xL}\alpha_{2p} (\boldsymbol{E_{0L}S_{xL}u})^T (I_{xL}\otimes H_{yL})  (\boldsymbol{E_{0L}S_{xL}u}) \\
&- (\boldsymbol{E_{0L}S_{xL}u})^T (I_{xL}\otimes H_{yL}) (\boldsymbol{u}-(E_{LR}\otimes I_{F2C})\boldsymbol{v}) \\
&+ \frac{\tau}{2}(\boldsymbol{u}-(E_{LR}\otimes I_{F2C})\boldsymbol{v})^T (I_{xL}\otimes H_{yL}) (\boldsymbol{u}-(E_{LR}\otimes I_{F2C})\boldsymbol{v}),\\
\boldsymbol{Q_3} =& h_{xR}\alpha_{2p} (\boldsymbol{E_{0R}S_{xR}v})^T (I_{xR}\otimes H_{yR}) (\boldsymbol{E_{0R}S_{xR}v}) \\
&- (\boldsymbol{E_{0R}S_{xR}v})^T (I_{xR}\otimes H_{yR}) ((E_{RL}\otimes I_{C2F})\boldsymbol{u}-\boldsymbol{v}) \\
&+ \frac{\tau}{2}((E_{RL}\otimes I_{C2F})\boldsymbol{u}-\boldsymbol{v})^T (I_{xR}\otimes H_{yR}) ((E_{RL}\otimes I_{C2F})\boldsymbol{u}-\boldsymbol{v}).
\end{split}
\end{equation*}
Since $\tilde{M}_{xL}$ and $\tilde{M}_{xR}$ are positive semidefinite, we have $\boldsymbol{Q_1}\geq 0$. To ensure $\boldsymbol{Q_2}\geq 0$ and $\boldsymbol{Q_3}\geq 0$, we need
\begin{equation*}
\begin{split}
& 2\sqrt{h_{xL}\alpha_{2p}}\sqrt{\tau/2}\geq 1 \text{ and } 2\sqrt{h_{xR}\alpha_{2p}}\sqrt{\tau/2}\geq 1. \\
\Rightarrow\quad &\tau\geq\max\left(\frac{1}{2\alpha_{2p}h_{xL}}, \frac{1}{2\alpha_{2p}h_{xR}}\right).
\end{split}
\end{equation*}
In the mesh shown in Figure \ref{Mesh_2_MC}, the mesh refinement ratio across the grid interface is $1:2$, i.e. $h_{xR}=\frac{1}{2}h_{xL}$. Thus,
\begin{equation}\label{tau}
\tau\geq\frac{1}{\alpha_{2p}h_{xL}}
\end{equation}
is the condition for $\boldsymbol{E_H^W}$ to be a discrete energy. Therefore, an energy estimate is obtained if (\ref{equality}), (\ref{inequality}) and (\ref{tau}) hold.

As shown above, in order to apply the energy method ($\ref{inequality}$) must be satisfied, which means that the energy estimate is only valid for the second and fourth order accurate schemes. The same interpolation operators are used for the Schr\"{o}dinger equation in \cite{Nissen2012} but ($\ref{inequality}$) is not needed for stability. For the Euler equations, ($\ref{inequality}$) is not needed when the standard SBP operators are used for the discretization, but needed when upwind SBP operators are used \cite{Mattsson2010}.  

\subsubsection{Convergence rate}
We discuss the accuracy properties of (\ref{Wave2dsemi}) by analyzing the truncation error of (\ref{Wave2dsemiu}). Equation (\ref{Wave2dsemiv}) can be analyzed in a similar way. The approximation error of the SBP operator $\boldsymbol{D_{2L}}$ is $\mathcal{O}(h^{2p})$ in the interior and $\mathcal{O}(h^p)$ near the interface, with the latter one being the dominant source of error. In the first penalty term $\boldsymbol{SAT_{u1}}$, a large interpolation error $\mathcal{O}(h^p)$ is located near the edge of the grid interface. Due to the $h^{-1}$ factor in both $\boldsymbol{H_{x}^{-1}}$ and $\boldsymbol{S_{xL}^T}$, the localized truncation error of $\boldsymbol{SAT_{u1}}$ is $\mathcal{O}(h^{p-2})$. Similarly, we find that the localized truncation errors of $\boldsymbol{SAT_{u2}}$ and $\boldsymbol{SAT_{\partial u}}$ are $\mathcal{O}(h^{p-2})$ and $\mathcal{O}(h^{p-1})$, respectively. Therefore, the localized truncation error of the semidiscretization (\ref{Wave2dsemiu}) is $\mathcal{O}(h^{p-2})$.

In \cite{Wang2015}, the convergence of the SBP--SAT discretization of the second order wave equation in one space dimension with a grid interface is analyzed. The result is that if the penalty parameter is chosen strictly larger than the limit value required for stability, the localized truncation error $\mathcal{O}(h^p)$ near the grid interface results in an error $\mathcal{O}(h^{p+1})$ in the solution for $p=1$, and an error $\mathcal{O}(h^{p+2})$ for $p\geq 2$. In other words, we gain one order in convergence if $p=1$ and two orders if $p\geq 2$.

In our case, the spatial dimension is two and there is the possibility of another gain in convergence. That is, the number of grid points with truncation error $\mathcal{O}(h^{p-2})$ is finite and independent of $h$. Hence, the L$_2$ norm of this truncation error is $\mathcal{O}(h^{p-1})$, and is one order higher than the pointwise truncation error. Therefore, we can hope to get an extra gain in convergence comparing with the corresponding one dimensional case. 

By a convergence test in \S\ref{sec-NE}, we find that  the extra gain is one order, which gives a total gain of two orders for $p=1$ and three orders for $p\geq 2$. That is, the localized truncation error $\mathcal{O}(h^{p-2})$ results in an error $\mathcal{O}(h^{p})$ in the solution for $p=1$, and an error $\mathcal{O}(h^{p+1})$ for $p\geq 2$. To obtain this convergence rate, it is important to choose the penalty parameter strictly larger than the value required for stability.

Comparing with the case of conforming grid interfaces, the convergence rate is one order lower. Even though a non--conforming grid interface allows for a local mesh refinement, the loss of accuracy may attenuate its efficiency in practice. To overcome the accuracy reduction by the non--conforming grid interfaces, we have tried to build interpolation operators with  error $\mathcal{O}(h^{2p})$ in the interior and error $\mathcal{O}(h^{p+1})$ near the edge of the interface, based on both diagonal and block norm SBP operators by using the symbolic software MAPLE and the approach presented in \cite{Mattsson2010}. However, we could not find a solution to the resulting system of equations. 

\subsection{An extension to T--junction interfaces}
\begin{figure}
\begin{center}
\subfloat[]{\includegraphics[width=4cm]{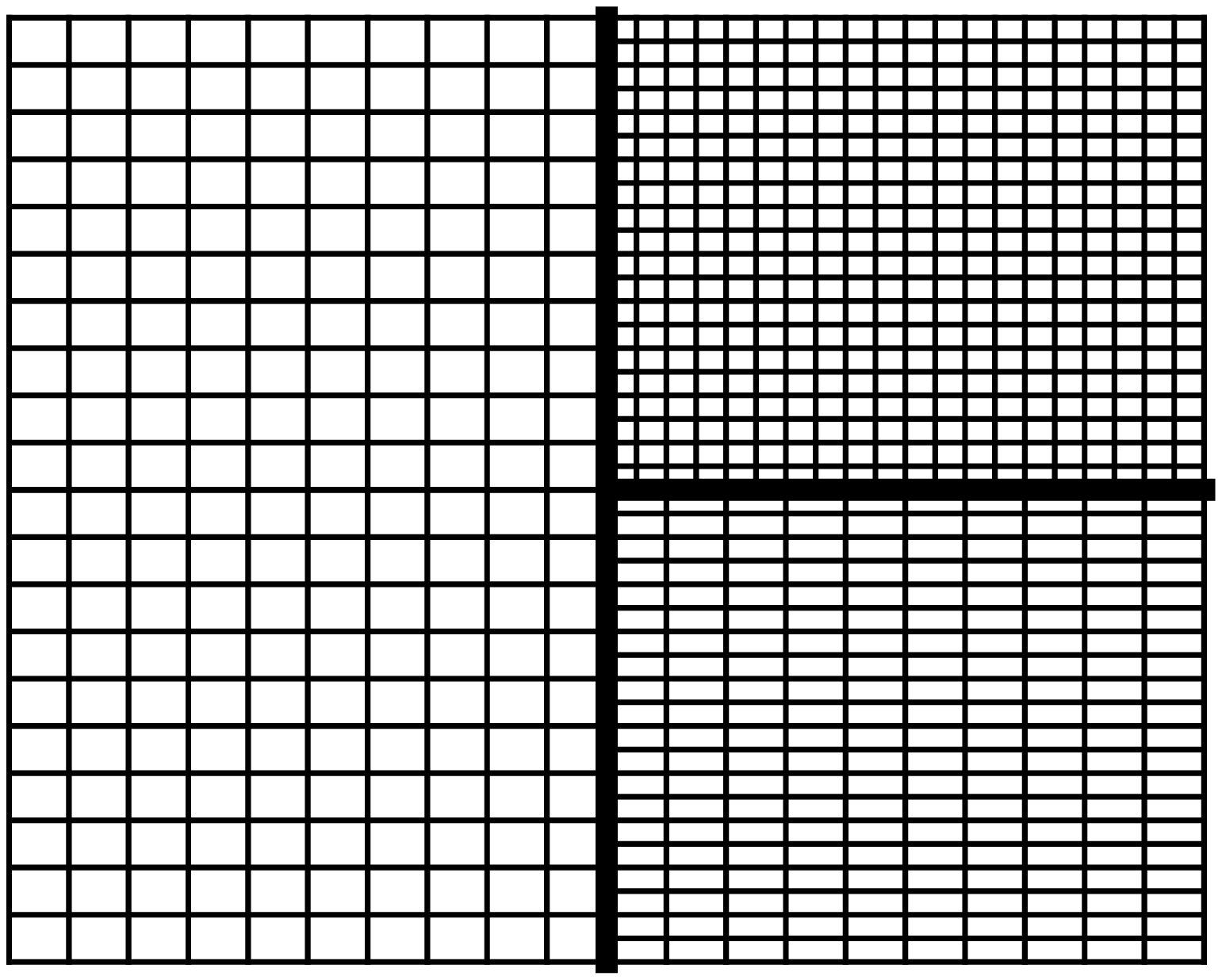}\label{Mesh_3_MC}}
\subfloat[]{\includegraphics[width=4cm]{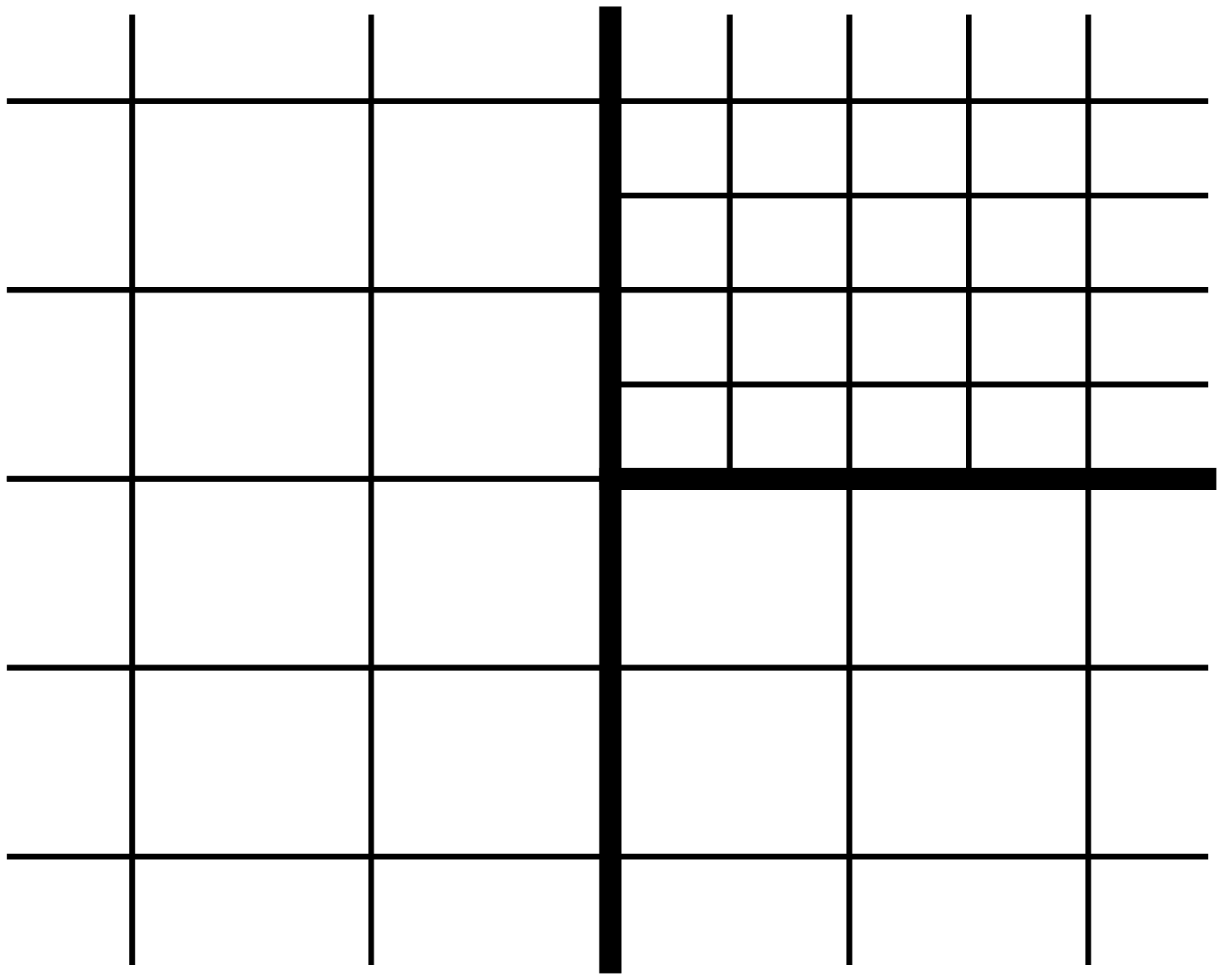}\label{Mesh_T1}}
\subfloat[]{\includegraphics[width=4cm]{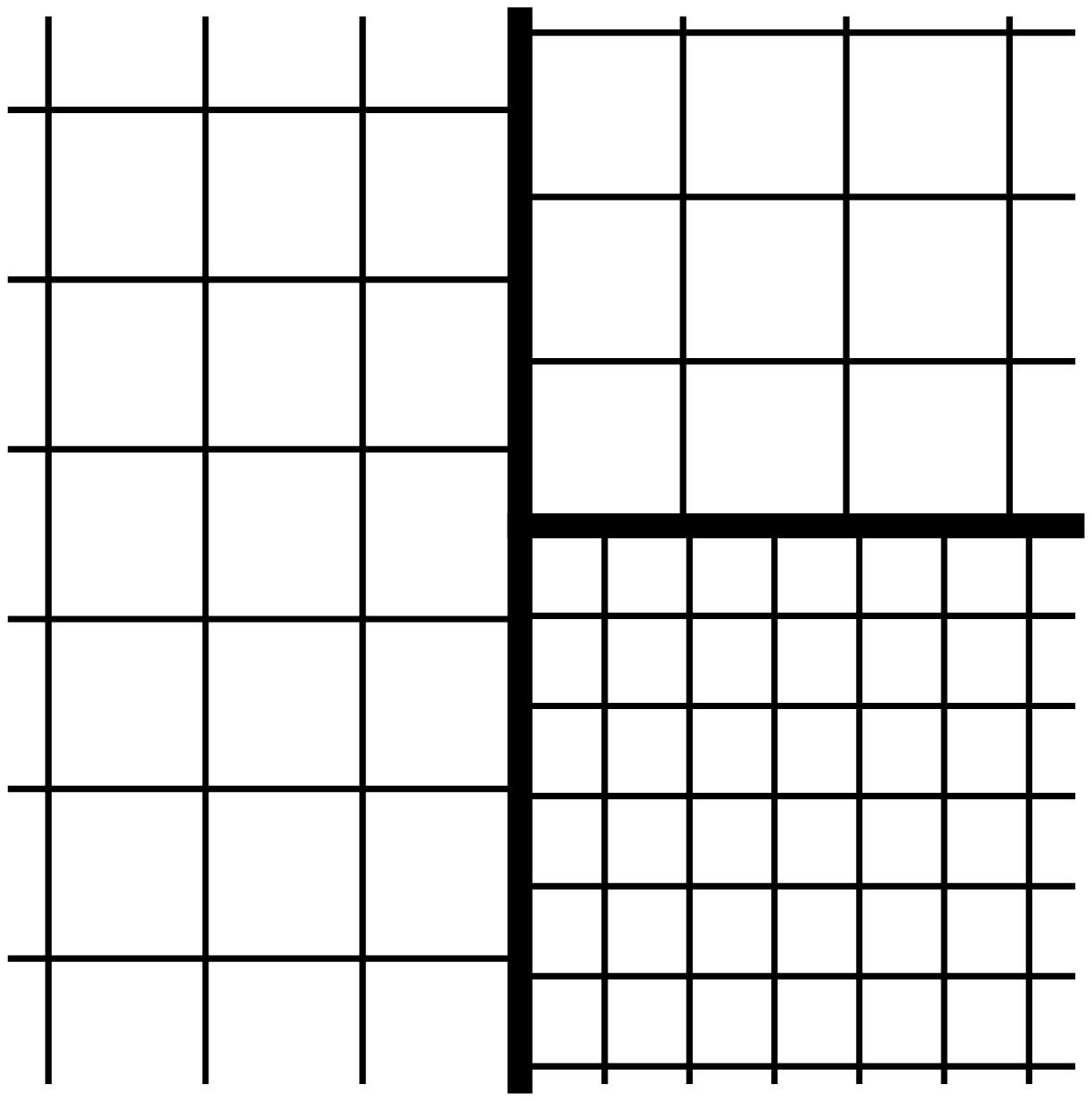}\label{NE_T}}
\caption{T--junction interfaces}
\label{Mesh_KW}
\end{center}
\end{figure}
In Figure \ref{Mesh_4_MC}, the interface between the two blocks on the left is conforming. It is then in many cases desirable not to consider it as a grid interface, but to use the mesh shown in Figure \ref{Mesh_3_MC} instead, where the interface there forms a T--junction. This is because SBP operators have a larger approximation error near the grid interface than that in the interior. The usage of redundant grid interfaces results in additional errors in the solution. Moreover, avoiding T--junction interfaces may end up in a bad partitioning of the computational domain with many unnecessary mesh blocks. With a straightforward application of the interpolation operators, instability occurs around the T--junction interface intersection point. The reason is that near the interface intersection point in Figure \ref{Mesh_3_MC}, the SBP norm in the vertical direction has the interior weights in the left domain, and the boundary weights in the two domains on the right. As a consequence, the norm--compatibility condition is violated and no energy estimate can be derived. 

In \cite{Nissen2015}, the T--junction operators are constructed to handle T--junction interfaces and are applied to the advection equation and the Schr\"{o}dinger equation in the SBP--SAT framework. Stability is proved by the energy method, but it comes with the cost that the T--junction operators introduce an error $\mathcal{O}(h^{2p})$ in the interior of the interface and $\mathcal{O}(h^{p})$ near the edge of the interface. The T--junction operators can also be used together with the interpolation operators to handle non--conforming grid interfaces. One constraint for the T--junction operators is that the interface intersection point must be a grid point in all involved mesh blocks, for example a close-up T--junction interface in Figure \ref{Mesh_T1}.  It is not straightforward to handle the T--junction interface shown in Figure \ref{NE_T} by the same technique.

\section{Non--conforming grid interfaces handled by projection operators}\label{sec-KW}
In \cite{Kozdon2015}, a new methodology of handling grid interfaces is introduced. In contrast to the interpolation operators which are based on a direct interpolation technique, the new methodology is based on a projection method. The highlights are that there is no strict requirement on the mesh refinement ratio, and the mesh blocks do not need to be conforming.

To illustrate how the projection method works, we consider again the mesh $\Omega$ shown in Figure \ref{Mesh_2_MC}, and denote $y^L$ in $\Omega_L$ and $y^R$ in $\Omega_R$ the grids on the interface. In addition, $z^L$ and $z^R$ denote the discrete finite difference solutions on $y^L$ and $y^R$. In a general setting, eight projection operators are used to move $z^L$ and $z^R$ between the two grids on the interface. Firstly, the discrete finite difference solution $z^L$ is projected by a projection operator $P_{f2p}^L$ to a piecewise continuous function based on the grid $y^L$. The associated mass matrix $M^L$ on $y^L$ is diagonal positive definite since Jacobi polynomials are used as the basis functions. Next, the glue grid $y$ with the mass matrix $M^g$ is defined as the grid that consists of the grid points on both $y^L$ and $y^R$. The projection operator $P_{p2g}^L$ is used to project the piecewise continuous function from $y^L$ to $y$, and is viewed as a basis transformation between polynomial spaces. Similarly, $P_{p2f}^L$ and $P_{g2p}^L$ are projection operators in the reversed direction corresponding to $P_{f2p}^L$ and $P_{p2g}^L$, respectively. $P_{f2p}^R$, $P_{p2f}^R$, $P_{p2g}^R$ and $P_{g2p}^R$ are the corresponding projection operators for the grid $y^R$.

Similar to the interpolation operators, there are norm--compatibility conditions for the projection operators:  
\begin{equation}\label{equ_pro}
H_{yL} P_{p2f}^L=(M^L P_{f2p}^L)^T,\quad H_{yR} P_{p2f}^R=(M^R P_{f2p}^R)^T,
\end{equation}
and
\begin{equation}\label{equ_pro_2}
M^g P_{p2g}^L=(M^L P_{g2p}^L)^T, \quad M^g P_{p2f}^R=(M^R P_{g2p}^R)^T.
\end{equation}
Define
\begin{equation}\label{I_Kozdon}
I_{L2R}^p=P_{p2f}^R P_{g2p}^R P_{p2g}^L P_{f2p}^L \text{ and } I_{R2L}^p=P_{p2f}^L P_{g2p}^L P_{p2g}^R P_{f2p}^R, 
\end{equation}
the operator $I_{L2R}^p$ moves $z^L$ from $y^L$ to $y^R$, and $I_{R2L}^p$ moves $z^R$ from $y^R$ to $y^L$. 

For the projection operators $P_{f2p}^{L/R}$ and $P_{p2f}^{L/R}$ which move the discrete finite difference solution to the subspace of piecewise continuous functions and back, the projection error is $\mathcal{O}(h^{2p})$ in the interior and $\mathcal{O}(h^p)$ near the edge, where $p=1,2,3,4,5$. $P_{p2g}^{L/R}$ and $P_{g2p}^{L/R}$ can be viewed as basis transformation operators between grids.  As a consequence, the projection error of $I_{L2R}^p$ and $I_{L2R}^p$ is also $\mathcal{O}(h^{2p})$ in the interior and $\mathcal{O}(h^p)$ near the edge, where $p=1,2,3,4,5$. In other words, $I_{L2R}^p$ and $I_{L2R}^p$ have the same accuracy properties as the diagonal norm SBP operators and the interpolation operators in \cite{Mattsson2010}.

\subsection{Interface treatment with the projection operators}
If no T--junction interface is present in the mesh, for example in Figure \ref{Mesh_2_MC} and \ref{Mesh_4_MC}, the operators $I_{L2R}^p$ and $I_{R2L}^p$ in (\ref{I_Kozdon}) are used to impose grid interface conditions in the same way as the interpolation operators discussed in \S\ref{sec-MC}. $I_{L2R}^p$ and $I_{R2L}^p$ satisfy an analogue of relation (\ref{equality}) up to tenth order accuracy, and relation (\ref{inequality}) up to fourth order accuracy. The stability analysis and accuracy properties of the interpolation operators in \S\ref{sec-MC} are still valid for $I_{L2R}^p$ and $I_{R2L}^p$. 

If a T--junction interface is present in the mesh, for example in Figure \ref{Mesh_KW}, we cannot use $I_{L2R}^p$ and $I_{R2L}^p$ in a direct way. Instability occurs around the junction point if on one side the SBP norm has the interior weights while on the other side it has the boundary weights, and no energy estimate can be obtained. To overcome the instability, the coupling is done on the glue grid. That is, we project finite difference solutions to the glue grid, compute the penalty terms there, and project them back to the finite difference grids. In this way, we avoid the instability caused by the SBP norms since the penalty terms are computed on the common glue grid.

The projection technique is very flexible to handle grid interfaces in the sense that we are free to choose the interface structure, the mesh refinement ratio and the accuracy of the diagonal norm SBP operators. In \cite{Kozdon2015}, the authors also couple the SBP finite difference method with the discontinuous Galerkin method, inspired by the relation between the discontinuous Galerkin spectral element method and the SBP--SAT finite difference method \cite{Gassner2013}. 

Finally, we remark that even with the mesh in Figure \ref{Mesh_2_MC} where both interpolation technique and projection technique are applicable, the interpolation operators $I_{C2F}(I_{F2C})$ are not the same as $I^p_{L2R}(I^p_{R2L})$ in (\ref{I_Kozdon}). The latter one has a wider stencil. In the construction procedure of these operators, one gets a system of linear equations after imposing stability and accuracy requirements. The solution of the linear system has a few free parameters. There are different ways to tune those free parameters. In \cite{Mattsson2010}, the free parameters are used to minimize the coefficients of the leading interpolation error in L$_2$ norm, while in \cite{Kozdon2015} the free parameters are used to minimize the distance between nearest eigenvalues of $P_{p2f}P_{f2p}$ for a finite difference grid of size 64. The choice of tuning free parameters has no influence on the theoretical order of accuracy, but may have an impact on condition (\ref{inequality}) and the practical accuracy. This is studied in more detail in the numerical experiments in \S\ref{sec-NE}.

\section{Numerical experiments}\label{sec-NE}
In this section, numerical experiments are performed to compare the schemes with the interpolation operators and the projection operators, and verify their stability and accuracy properties. Moreover, we also conduct two numerical experiments to study the efficiency of local mesh refinement by solving the wave equation on a domain with a complex geometry. 

The L$_2$ error and maximum error are computed as the norm of the difference between the exact solution $u^{ex}$ and the numerical solution $u^h$ according to
\begin{equation*}
\begin{split}
&\|u^h-u^{ex}\|_{\text{L}_2}=\sqrt{h^d(u^h-u^{ex})^T(u^h-u^{ex})}, \\
&\|u^h-u^{ex}\|_{\infty}=\max |u^h-u^{ex}|/\text{amp}, 
\end{split}
\end{equation*}
where $d$ is the dimension of the equation and amp is the maximum amplitude of the solution. The convergence rate is computed by
\begin{equation*}
q=\log\left(\frac{\|u^h-u^{ex}\|}{\|u^{2h}-u^{ex}\|}\right) \bigg/ \log\left(\frac{1}{2}\right).
\end{equation*}

\subsection{Stability study}\label{sec-stability}
We begin with an eigenvalue analysis for condition (\ref{inequality}). The computational domain is $x\in [-1,1]$ and $y\in [0,1]$ with a grid interface at $x=0$. In the left domain the number of grid points is 26 in both $x$ and $y$ directions, while in the right domain the number of grid points is 51 in both $x$ and $y$ directions. The mesh refinement ratio is $1:2$ across the grid interface, which enables both the interpolation operators and projection operators applicable. The matrices $\Xi_L$ and $\Xi_R$ are symmetric, so they have real eigenvalues. We denote $k_L$ and $k_R$ the smallest eigenvalue of $\Xi_L$ and $\Xi_R$ in (\ref{inequality}), scaled by the mesh size:
\begin{equation*}
k_L=\min(\text{eig}(\Xi_L))/h_{yL} \quad k_R=\min(\text{eig}(\Xi_R))/h_{yR}.
\end{equation*}
The reason for the scaling is that the elements in $\Xi_L / \Xi_R$ are proportional to $h_{yL} / h_{yR}$.

In Table \ref{StabilityNE}, we list $k_L$ and $k_R$ for the interpolation operators and the projection operators in Column three and four, respectively. For the interpolation operators, $(\ref{inequality})$ holds for both the second and fourth order accurate cases with errors up to machine precision. For the second order accurate case, we can prove  that $k_L,k_R\geq 0$ independent of $h$, because $\Xi_L$ is diagonally dominant and $\Xi_R$ can be transformed to a diagonally dominant matrix without changing the signs of the eigenvalues. For the sixth and eighth order cases, $(\ref{inequality})$ no longer holds. The difference between these two cases is that $k_L$ and $k_R$ are close to zero for the sixth order case, but far away from zero with the eighth order case. When increasing the number of grid points, the values of $k_L$ and $k_R$ remain unchanged.

For the projection operators, $(\ref{inequality})$ also holds also for the second and fourth order accurate cases. For the sixth, eighth and tenth order cases, $(\ref{inequality})$ does not hold anymore but the values of $k_L$ and $k_R$ are close to zero, and they become slightly closer to zero as the mesh is refined.  
\begin{table}
\centering
\begin{tabular}{c c c c c c}
\toprule
$2p$                    &            & IO                              & PO                            & IO                             & PO \\ \midrule
\multirow{2}{*}{2} & $k_L$ & $-6.2\cdot 10^{-17}$ & $-7.7\cdot 10^{-16}$ & \multirow{2}{*}{$-5.0\cdot 10^{-3}$} & \multirow{2}{*}{$-5.0\cdot 10^{-3}$} \\ \cline{2-4}
                            & $k_R$ & $-1.0\cdot 10^{-16}$ & $-8.5\cdot 10^{-16}$ & & \\ \midrule
\multirow{2}{*}{4} & $k_L$ & $-6.7\cdot 10^{-17}$ & $-8.3\cdot 10^{-16}$ & \multirow{2}{*}{$-5.0\cdot 10^{-3}$} & \multirow{2}{*}{$-5.0\cdot 10^{-3}$} \\ \cline{2-4}
                            & $k_R$ & $-1.4\cdot 10^{-16}$ & $-6.3\cdot 10^{-16}$ & & \\ \midrule
\multirow{2}{*}{6} & $k_L$ & $-6.9\cdot 10^{-1}$ & $-7.5\cdot 10^{-5}$ & \multirow{2}{*}{$-5.0\cdot 10^{-3}$} & \multirow{2}{*}{$-5.0\cdot 10^{-3}$} \\ \cline{2-4}
                            & $k_R$ & $-8.0\cdot 10^{-1}$ & $-8.6\cdot 10^{-5}$ & & \\ \midrule
\multirow{2}{*}{8} & $k_L$ & $-3.1\cdot 10^{1}$ & $-4.4\cdot 10^{-4}$ & \multirow{2}{*}{$8.7\cdot 10^{3}$} & \multirow{2}{*}{$-5.0\cdot 10^{-3}$} \\ \cline{2-4}
                            & $k_R$ & $-7.4\cdot 10^{1}$ & $-4.8\cdot 10^{-4}$ & & \\ \midrule
\multirow{2}{*}{10} & $k_L$ & ------ & $-4.8\cdot 10^{-4}$ & \multirow{2}{*}{------} & \multirow{2}{*}{$-5.0\cdot 10^{-3}$} \\ \cline{2-4}
                            & $k_R$ & ------ & $-1.1\cdot 10^{-3}$ & & \\ \bottomrule
\end{tabular}
\caption{IO: interpolation operators, PO: projection operators. Column 3 and 4 correspond to the numerical study of relation (\ref{inequality}), and Column 5 and 6 correspond to the eigenvalue analysis of (\ref{ode}).}
\label{StabilityNE}
\end{table}

Another way of analyzing stability through numerical experiments is to write the semidiscretized equation (\ref{Wave2dsemi}) as a system of ordinary differential equations 
\begin{equation}\label{ode}
\boldsymbol{z_{tt}=Qz+F}.
\end{equation}
It is stable if the eigenvalues of $\boldsymbol{Q}$ are real and non--positive. Otherwise, the numerical scheme is unstable. We have computed the eigenvalues of $\boldsymbol{Q}$ by using the same mesh as for the computation of $k_L$ and $k_R$. All the eigenvalues are real. In Table \ref{StabilityNE}, the largest eigenvalue of $\boldsymbol{Q}$ is shown in Column five and six for the schemes with the interpolation operators and projection operators, respectively. For the numerical scheme with the interpolation operators, it is stable for the second, fourth and sixth order cases, and unstable for the eighth order case. For the numerical scheme with the projection operators, it is stable for up to tenth order cases even though relation (\ref{inequality}) only holds for up to fourth order scheme. The eigenvalue analysis indicates that (\ref{inequality}) is sufficient but not necessary for stability.

\subsection{Accuracy study}
In this section, the convergence of the SBP--SAT method applied to the wave equation (\ref{Wave2d}) with a non--conforming grid interface is investigated. The analytical solution to (\ref{Wave2d}) is manufactured, which means that a closed form is chosen and is used to obtain the initial and boundary data. At the outer boundaries, Dirichlet boundary condition is imposed weakly by the SAT method as described in \cite{Mattsson2009}.

\subsubsection{A non--conforming grid interface}
In the first numerical experiment, the computational domain is $[-1,1]\times [0,1]$ where a grid interface is located at $x=0, y\in [0,1]$. In the grid refinement level $r=0$, the numbers of grid points in the left block and right block are $26\times 26$ and $51\times 51$. The mesh sizes are halved in both $x$ and $y$ directions when $r$ is increased by one. In this setting, the grid refinement ratio is $1:2$, and both the interpolation operators and the projection operators are applicable for the numerical treatment of interface conditions. 
The fourth order Runge--Kutta method is used as the time integrator with the step size in time chosen so small that the temporal error is negligible compared with the spatial error. 

The manufactured solution to (\ref{Wave2d}) is chosen to be
\begin{equation}\label{Analytical}
U(x,y,t)=\cos(5x+1)\cos(5y+2)\cos(5\sqrt{2}t+3).
\end{equation}
 
The computational results are shown in Table \ref{MC_table}, where $2p$ and $r$ in the first two columns denote the order of accuracy and the mesh refinement level, respectively. In Column 3, 4 and 5 the errors in L$_2$ norm, the convergence rates in L$_2$ norm and maximum norm are shown for the numerical schemes with the interpolation operators, whereas in Column 6, 7 and 8 the corresponding results obtained by the schemes with the projection operators are shown.
\begin{table}
\centering
\begin{tabular}{c c c c c c c c}
\toprule
$2p$ & $r$ & L$_2$ error & $q_{L_2}$ & $q_M$ & L$_2$ error & $q_{L_2}$ & $q_M$ \\ \midrule
2 &  0   & $1.28\cdot 10^{-2}$ & ------ & ------  & $3.76\cdot 10^{-2}$ & ------ & ------   \\ \midrule
   &  1   & $7.05\cdot 10^{-3}$ & 0.87 &  0.85 & $2.02\cdot 10^{-2}$ & 0.90  & 0.83 \\ \midrule
   &  2   & $3.70\cdot 10^{-3}$ & 0.93 &  0.89 & $1.04\cdot 10^{-2}$ & 0.95  &0.91 \\  \midrule
   &  3   & $1.90\cdot 10^{-3}$ & 0.96 &  0.95 & $5.31\cdot 10^{-3}$ & 0.98  &0.95 \\ \midrule
4 &  0   & $8.41\cdot 10^{-4}$ & ------ & ------ & $1.22\cdot 10^{-3}$ & ------ & ------ \\ \midrule
   &  1   & $1.18\cdot 10^{-4}$ & 2.83 & 2.57 & $1.15\cdot 10^{-4}$ & 3.40  &3.46\\ \midrule
   &  2   & $1.52\cdot 10^{-5}$ & 2.96 & 2.86 &  $1.31\cdot 10^{-5}$ & 3.14  &2.42\\  \midrule
   &  3   & $1.91\cdot 10^{-6}$ & 2.99 & 2.69 &  $1.61\cdot 10^{-6}$ & 3.03  &1.49\\ \midrule
6 &  0   & $7.65\cdot 10^{-5}$ & ------ & ------ & $1.10\cdot 10^{-4}$ & ------  & ------ \\ \midrule
   &  1   & $6.91\cdot 10^{-6}$ & 3.47 & 2.92 &  $8.94\cdot 10^{-6}$ & 3.62  &2.62\\ \midrule
   &  2   & $5.02\cdot 10^{-7}$ & 3.78 & 2.56 &  $5.88\cdot 10^{-7}$ & 3.93  &3.97\\  \midrule
   &  3   & $3.33\cdot 10^{-8}$ & 3.91 & 3.01 &  $3.85\cdot 10^{-8}$ & 3.93  &3.01\\ \midrule
8 &  0   &                        &  &  &  $7.21\cdot 10^{-4}$ & ------  & ------  \\ \midrule
   &  1   &                        &  &  &  $1.53\cdot 10^{-5}$ & 5.56 & 4.73\\ \midrule
   &  2   &                         &  &  & $3.97\cdot 10^{-7}$ & 5.27 & 4.16\\  \midrule
   &  3   &                       &  & &    $1.05\cdot 10^{-8}$ & 5.25 & 4.95\\ \midrule
10&  0   &                          & &  &  $3.39\cdot 10^{-5}$ & ------  & ------       \\ \midrule
   &  1   & & &  & $5.12\cdot 10^{-7}$ & 6.05  & 5.08  \\ \midrule
   &  2   & & &  & $1.19\cdot 10^{-8}$ & 5.43  & 4.54 \\  \midrule
   &  3   & & &  & $2.27\cdot 10^{-10}$ & 5.71  & 4.94  \\ \bottomrule
\end{tabular}
\caption{Convergence of the SBP--SAT scheme for the wave equation with a grid interface. The interface conditions are handled by the interpolation operators (Column 3,4,5) and the projection operators  (Column 6,7,8).}
\label{MC_table}
\end{table}

For the scheme with the interpolation operators and time step $\Delta t=0.1h$, the convergence rates in L$_2$ norm are 1, 3 and 4 for the second, fourth and sixth order schemes, respectively. This agrees with our accuracy discussion in \S\ref{sec-MC}. Though stability cannot be proved by the energy method for the sixth order accurate case, it seems that for this particular setting the scheme is stable and exhibits the expected convergence rate. Instability occurs when using the eighth order accurate scheme.  

With the projection operators, the convergence rate in L$_2$ norm is one for the second order accurate scheme, and $p+1$ for the fourth, sixth, eighth and tenth order accurate schemes, which agrees with the  accuracy discussion in \S\ref{sec-KW}. We note that though an energy estimate does not exist for the sixth, eighth and tenth order accurate cases, the schemes are stable and converge as expected. The time step is $\Delta t=0.1h$ for $2p=2,4,6$. With this time step, the tenth order accurate scheme yields slightly lower convergence rate than expected, and the result shown in Table \ref{MC_table} is obtained with  $\Delta t=0.05h$. The eighth order accurate scheme is a special one, since with $\Delta t=0.05h$ it is even unstable. To obtain the results in Table \ref{MC_table}, $\Delta t=0.025h$ is used as the time step, which indicates that the eighth order accurate semidiscretized equation is stiff. Moreover, the error obtained with the eighth order accurate scheme is larger than the error obtained with the sixth order accurate scheme, except for the finest mesh refinement level.  

From Table \ref{MC_table}, it is also observed that the errors are similar to each other for the schemes of the same order of accuracy with interpolation operators and projection operators. 

\subsubsection{A T--junction interface}\label{sec-T}
Next, we consider the computational domain $[-1,1]^2$ that is divided into three mesh blocks as shown in Figure \ref{Mesh_3_MC}. The interfaces are located at $x=0,y\in [-1,1]$ and $y=0,x\in [0,1]$. In the grid refinement level $r=0$, the numbers of grid points in Block 1 (left), Block 2 (right--up) and Block 3 (right--down) are $28\times 51$, $27\times 25$ and $51\times 50$, respectively. The mesh sizes are halved in both $x$ and $y$ directions when $r$ is increased by one. This partitioning and meshing result in a highly non--conforming grid interface with a close-up shown in Figure \ref{NE_T}. The interface conditions are imposed weakly by the SAT method with the projection operators. To test convergence, (\ref{Analytical}) is used as the analytical solution. The computational results are shown in table \ref{T_table}.

\begin{table}
\centering
\begin{tabular}{c c  c  c c c}
\toprule
$2p$ & $r$ & L$_2$ error & $q_{L_2}$ & Maximum error & $q_M$ \\ \midrule
2 &  0   & $1.07\cdot 10^{-1}$ & ------ & $1.40\cdot 10^{-1}$ &  ------ \\ \midrule
   &  1   & $5.83\cdot 10^{-2}$ & 0.87 & $8.21\cdot 10^{-2}$ & 0.77 \\ \midrule
   &  2   & $3.04\cdot 10^{-2}$ & 0.94 & $4.47\cdot 10^{-2}$ & 0.88 \\  \midrule
   &  3   & $1.55\cdot 10^{-2}$ & 0.97 & $2.35\cdot 10^{-2}$ & 0.93 \\ \midrule
4 &  0   & $3.03\cdot 10^{-3}$ & ------ & $9.39\cdot 10^{-3}$ &  ------ \\ \midrule
   &  1   & $3.08\cdot 10^{-4}$ & 3.29 & $2.13\cdot 10^{-3}$ & 2.14 \\ \midrule
   &  2   & $3.68\cdot 10^{-5}$ & 3.07 & $5.34\cdot 10^{-4}$ & 1.99 \\  \midrule
   &  3   & $4.37\cdot 10^{-6}$ & 3.07 & $1.33\cdot 10^{-4}$ & 2.01 \\ \midrule
6 &  0   & $3.00\cdot 10^{-4}$ & ------ & $2.04\cdot 10^{-3}$ &  ------ \\ \midrule
   &  1   & $1.63\cdot 10^{-5}$ & 4.20 & $2.23\cdot 10^{-4}$ & 3.19 \\ \midrule
   &  2   & $9.47\cdot 10^{-7}$ & 4.11 & $3.26\cdot 10^{-5}$ & 2.77 \\  \midrule
   &  3   & $6.22\cdot 10^{-8}$ & 3.93 & $4.49\cdot 10^{-6}$ & 2.86 \\ \midrule
8 &  0   & $1.04\cdot 10^{-2}$ & ------ & $1.46\cdot 10^{-2}$ &  ------ \\ \midrule
   &  1   & $6.06\cdot 10^{-5}$ & 7.43 & $3.52\cdot 10^{-4}$ & 5.38 \\ \midrule
   &  2   & $1.36\cdot 10^{-6}$ & 5.48 & $5.74\cdot 10^{-6}$ & 5.94 \\  \midrule
   &  3   & $3.83\cdot 10^{-8}$ & 5.15 & $3.92\cdot 10^{-7}$ & 3.87 \\ \midrule
10&  0   & $1.12\cdot 10^{-4}$ & ------ & $8.56\cdot 10^{-4}$ &  ------ \\ \midrule
   &  1   & $1.51\cdot 10^{-6}$ & 6.23 & $2.15\cdot 10^{-5}$ & 5.31 \\ \midrule
   &  2   & $2.31\cdot 10^{-8}$ & 6.03 & $1.14\cdot 10^{-6}$ & 4.24 \\  \midrule
   &  3   & $3.20\cdot 10^{-10}$ & 6.17 & $1.95\cdot 10^{-8}$ & 5.86 \\ \bottomrule
\end{tabular}
\caption{Convergence for the wave equation with a T--junction interface handled by the projection operators.}
\label{T_table}
\end{table}
Clearly, $(p+1)^{th}$ convergence rate in L$_2$ norm is obtained for $p=2,3,4,5$ and first order convergence rate is obtained for $p=1$. Here, we observe again that the schemes higher than fourth order accuracy are stable though no energy estimate can be obtained.

\subsection{Efficiency study}
\label{sec:efficiency_study}
In many applications, the frequencies of the present waves are given by initial and boundary data, and internal forcing functions. The wavelength of a wave is determined by the ratio between the wave speed of the material in which the wave is traveling and  the frequency of the wave. The accuracy of a numerical solution can be stated in terms of how many grid points per wavelength are used to resolve the present waves \cite{Kreiss1972}. A reduction in wave speed confined to a subset of the physical domain yields waves of a shorter wavelength localized in that subset. For accuracy it is therefore necessary to refine the grid to compensate for the shorter present wavelengths. For computational efficiency it is important that this refinement is done only in the subset that constitutes the slower media, since it is only in the slower media that wavelengths are reduced. As an example Figure \ref{fig:IS} shows the scattering and diffraction of an acoustic time--harmonic plane wave impinging on a circular region of a slower material, the wavelength is seen to be reduced inside the circular region.
\begin{figure}
        \includegraphics[trim=0cm 3.5cm 0cm 2.5cm, clip=true, width= \textwidth]{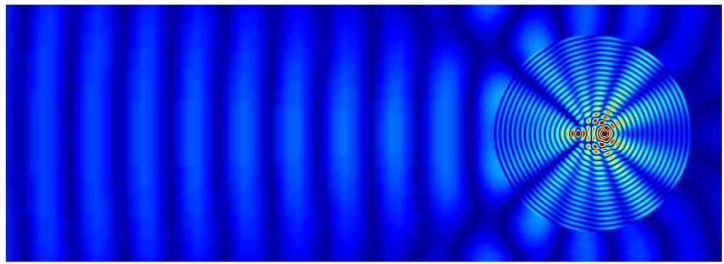} 
        	\caption{An example of an acoustic time--harmonic plane wave impinging on a circular inclusion, where the wavelength is much smaller inside the circular inclusion than that outside.}
	\label{fig:IS}
\end{figure}

Geometrical features of the physical domain such as a complex boundary or an internal cavity introduce a local radius of curvature. A small local radius of curvature compared with the present wavelengths can imply difficulties when generating a computational grid. As an example Figure \ref{fig:CS} shows the scattering of an acoustic time-harmonic plane wave impinging on a circular cavity of a radius of curvature smaller than the wave length of the incoming and scattered waves. A part of the grid used to represent the wave field is shown as an inset, where it is seen in Figure \ref{fig:MeshC} that the quality of the grid is impaired as the grid spacing gets unnecessarily small close to the cavity.       
\begin{figure}
	\subfloat[A circular cavity]{\includegraphics[width= 0.49\textwidth]{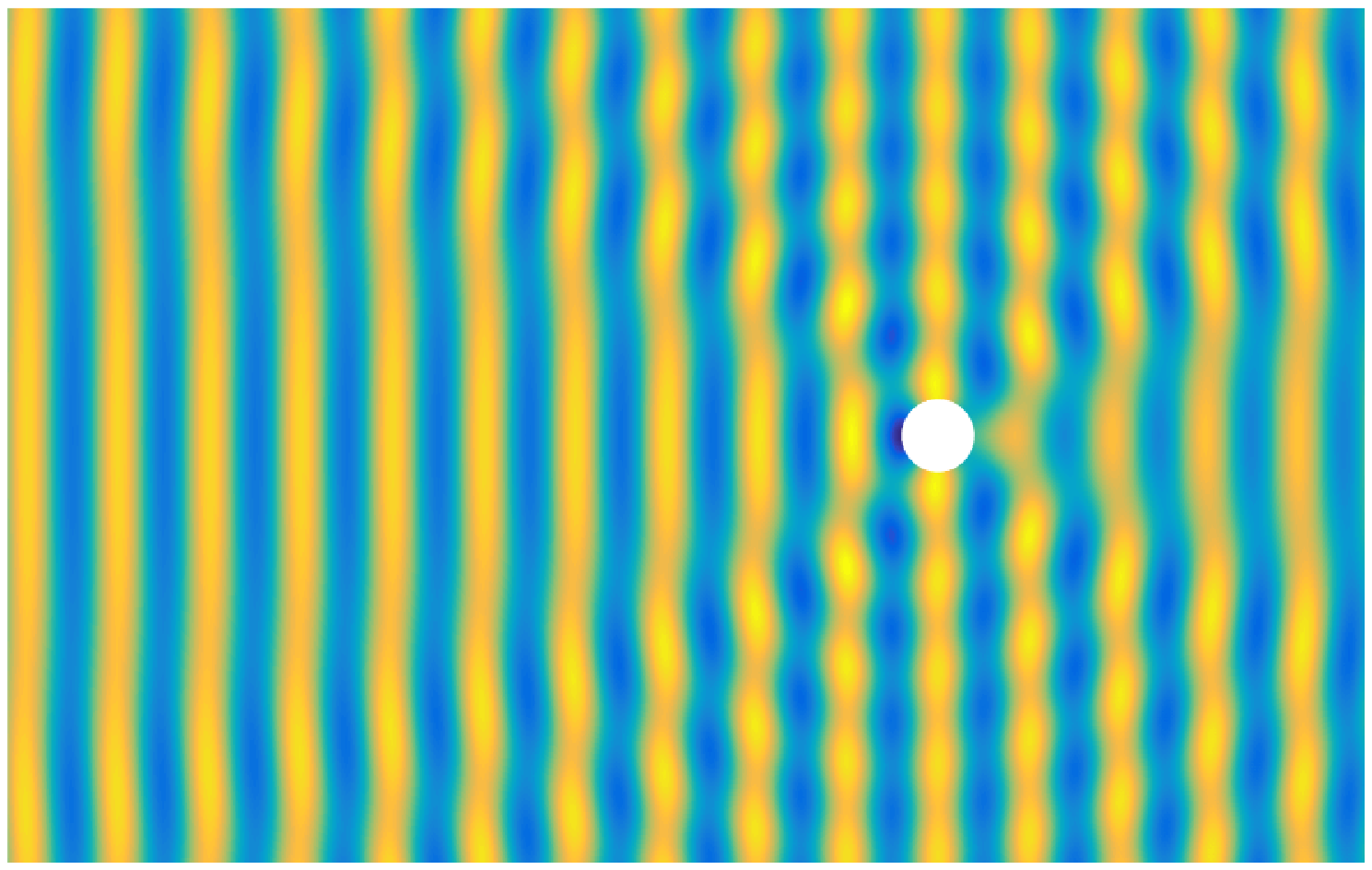}\label{fig:CS}}
	\subfloat[A close-up of the mesh near the cavity]{\includegraphics[width= 0.49\textwidth]{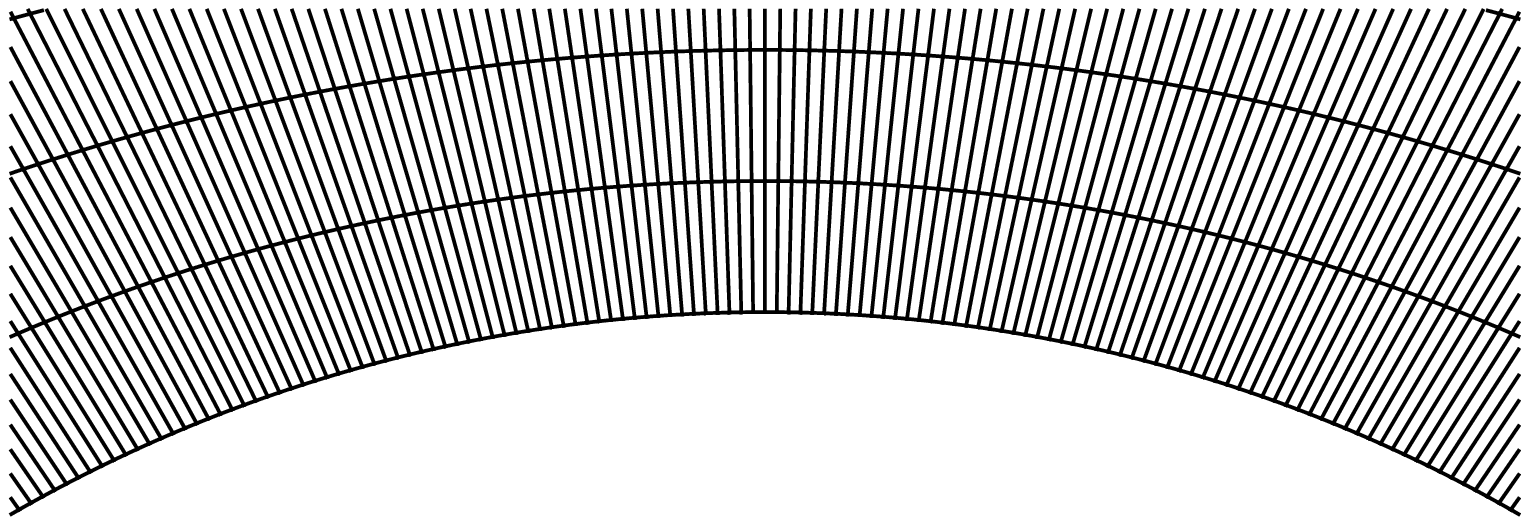}\label{fig:MeshC}}
	\caption{An examples of an acoustic time--harmonic plane wave impinging on a circular inclusion cavity with a grid close to the cavity where the grid size in the azimuth direction is much smaller than that in the radial direction.}
\end{figure}

In the preceding experiments it has been verified that using interpolation and projection operators to patch together the computational grid in a multi--block fashion yields a stable discretization, the convergence rate, however, was seen to be reduced. In the following experiments we investigate the practical benefit of using interpolation operators and projection operators in a region with a slower wave speed or a geometrical feature of a small radius of curvature albeit the order reduction. In particular, we will consider experiments involving acoustic waves impinging on a circular cavity and a circular inclusion of a differing material.

The numerical method used to solve the acoustic wave equation in the following experiments is based on the SBP-SAT scheme described in \cite{Virta2014}. The geometrical features are handled by using a multi-block strategy to decompose the physical domain into blocks, where each block allows for a mapping to curvilinear coordinates. In \cite{Virta2014} the blocks that constitute the domain are discretized by using conforming grids and patched together by the SAT method. In this paper we allow for non--conforming grids by implementing interpolation operators as well as projection operators into the handling of the multi-block interfaces.    

The following numerical experiments use two different two dimensional domains:
\begin{itemize}
	\item $\mathcal{D}_1$: An acoustic plane with a circular cavity of radius $a$.
	\item $\mathcal{D}_2$: An acoustic plane with a circular inclusion of radius $a$.
\end{itemize}
 The geometries of the domains are handled by decomposing each domain in a multi--block fashion. The blocks are then patched together to compound the entire domain. A detailed description of how these two grids are constructed is presented in the Appendix.
\subsubsection{A circular cavity}
In this numerical experiment, we consider a domain of an infinite homogeneous medium with a circular cavity of radius $a=1$. Let a plane harmonic wave $u_i=e^{i(\omega t-\gamma x)}$ propagate in that domain and impinge on the cavity. A scattered wave $u_s$ is generated when the incident wave hits the cavity, and the total displacement $u_i+u_s$ satisfies the wave equation. Homogeneous Neumann boundary condition is imposed at the cavity boundary. A detailed derivation of an analytical solution is found in \cite[\S 7]{Graff1991}. 

We take $\omega=2\pi$ and $c=2.5$, which give a wavelength 2.5. The computational domain $\mathcal{D}_1$ is chosen to be the rectangular $[-25.5,11.7]\times[-11.7,11.7]$ and the cavity is centred at the origin. Two ways of partitioning the domain are considered, namely the N--partitioning and the T--partitioning shown in Figure \ref{fig:naive} and \ref{fig:nonnaive}, respectively. In the N--partitioning approach, we only use conforming gird interfaces and conforming mesh blocks. The cavity is surrounded by four blocks that constitute the square $[-11.7,11.7]^2$,  which is attached by a rectangular domain to the left. The numbers of grid points on each edge are shown in the figure, and are chosen so that approximately 20 grid points per wavelength are used in the discretization. In this setting, the mesh is of bad quality since the mesh size near the cavity is significantly smaller than that near the outer boundaries. To overcome this drawback, we propose the T--partitioning where the cavity is surrounded by a small square block $[-1.3,1.3]^2$. Here, all the grid interfaces are also conforming but a T--junction interface is present at $x=-11.7$ with the intersection points marked by the dots. Again we choose the mesh size so that there are approximately 20 grid points per wavelength, and here it is only over--resolved in the small block $[-1.3,1.3]^2$. The T--partitioning results in a mesh of 54903 grid points.  The number of grid points with the N--partitioning is about doubled to 109867.

\begin{figure}
	\subfloat[N--partitioning]{\includegraphics[width=0.49\textwidth]{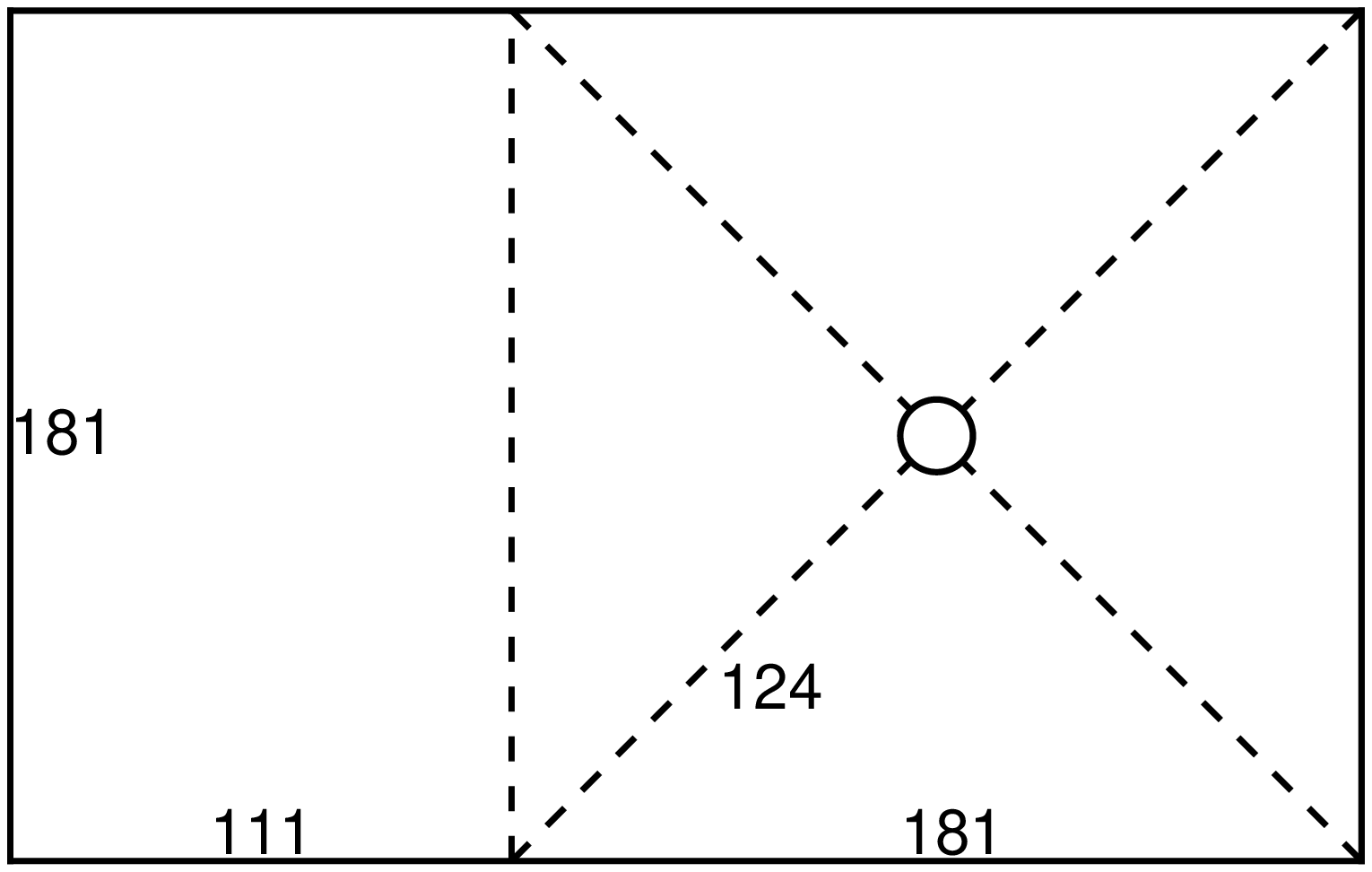}\label{fig:naive}}
	\subfloat[T--partitioning]{\includegraphics[width=0.49\textwidth]{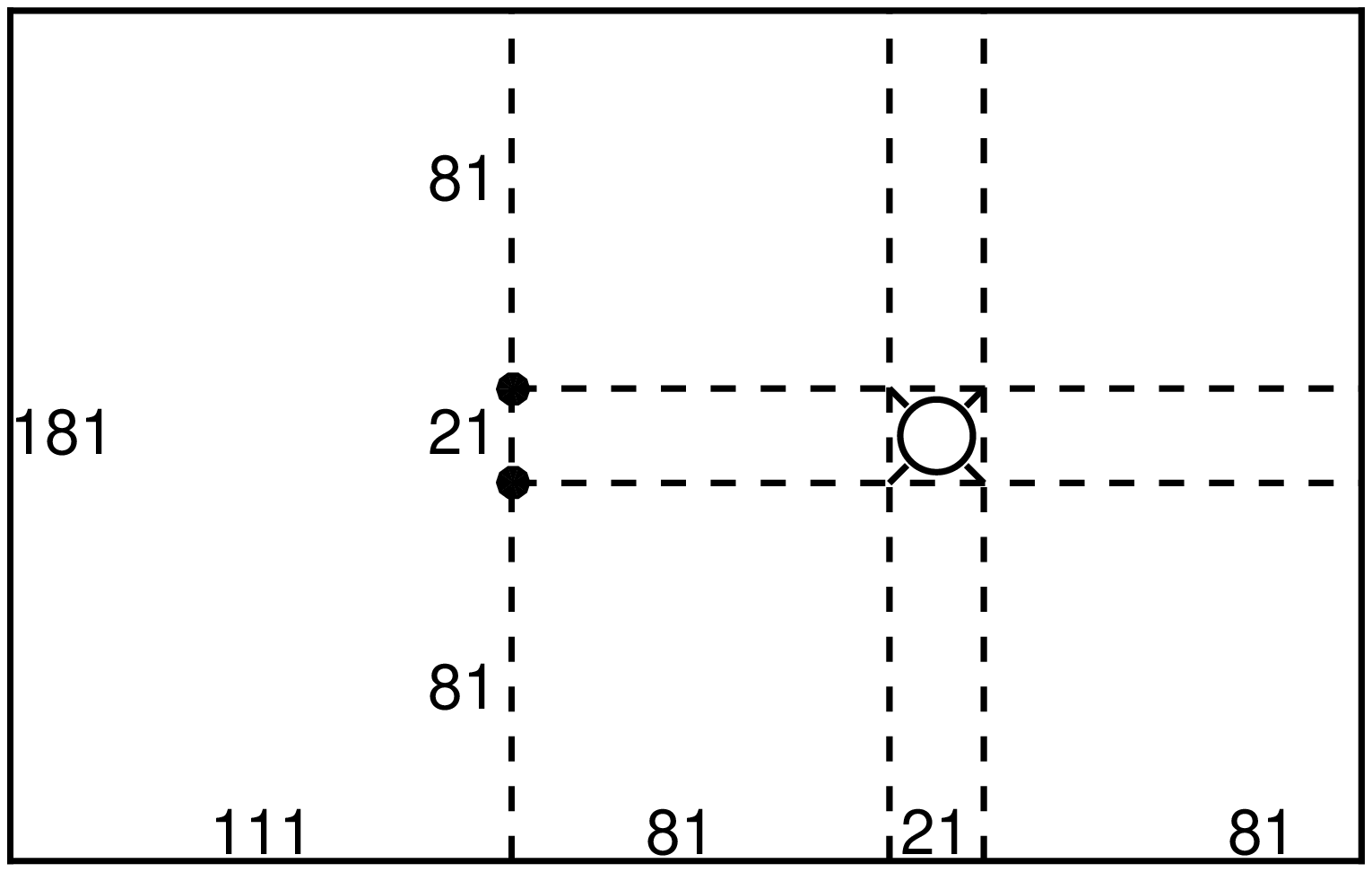}\label{fig:nonnaive}}
	\caption{Computational domain of the circular cavity experiment }
\end{figure}

\begin{figure}
	\subfloat[Fourth order accurate scheme]{\includegraphics[width=0.49\textwidth]{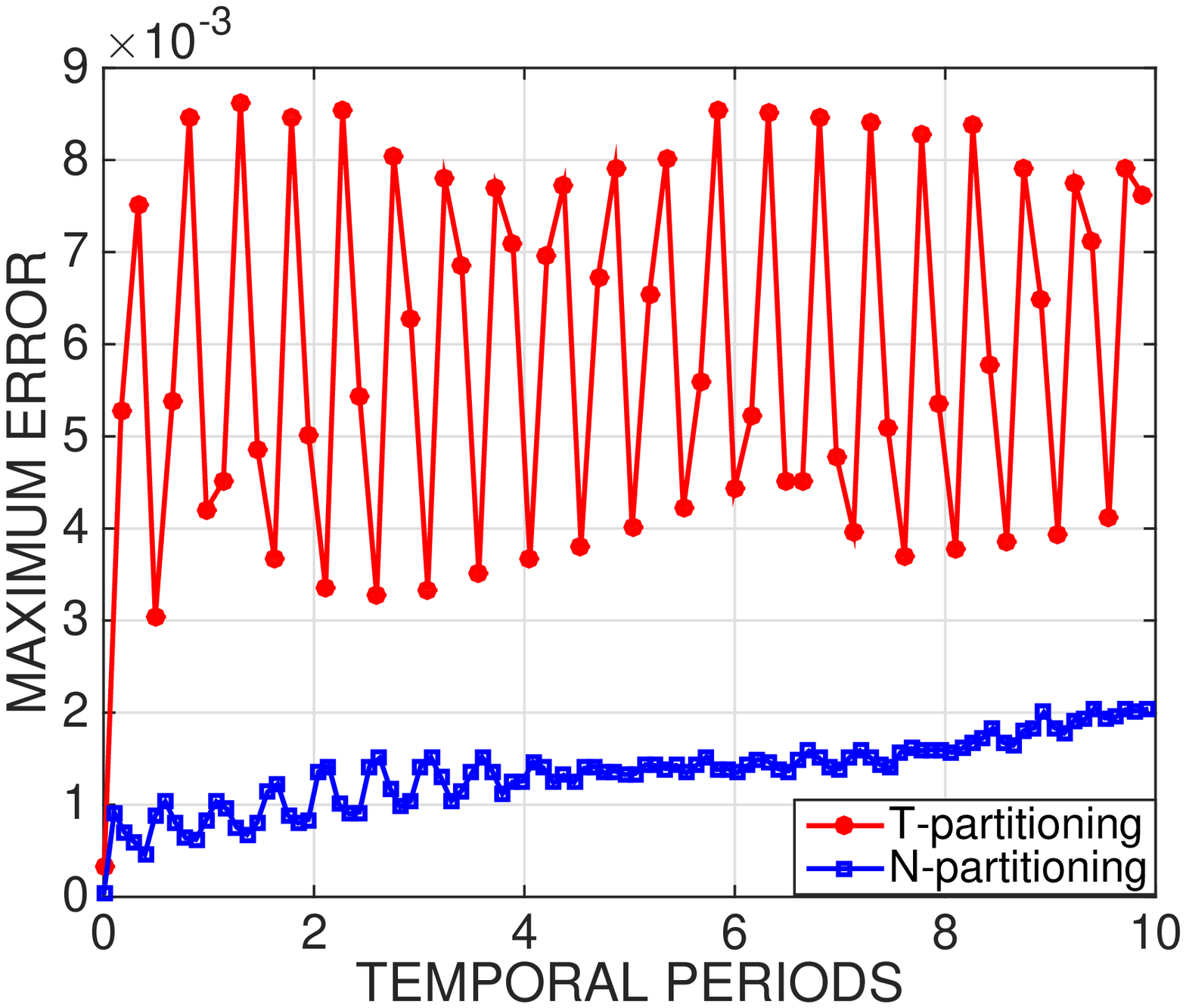}\label{fig:FourthCavity}}
	\subfloat[Sixth order accurate scheme]{\includegraphics[width=0.49\textwidth]{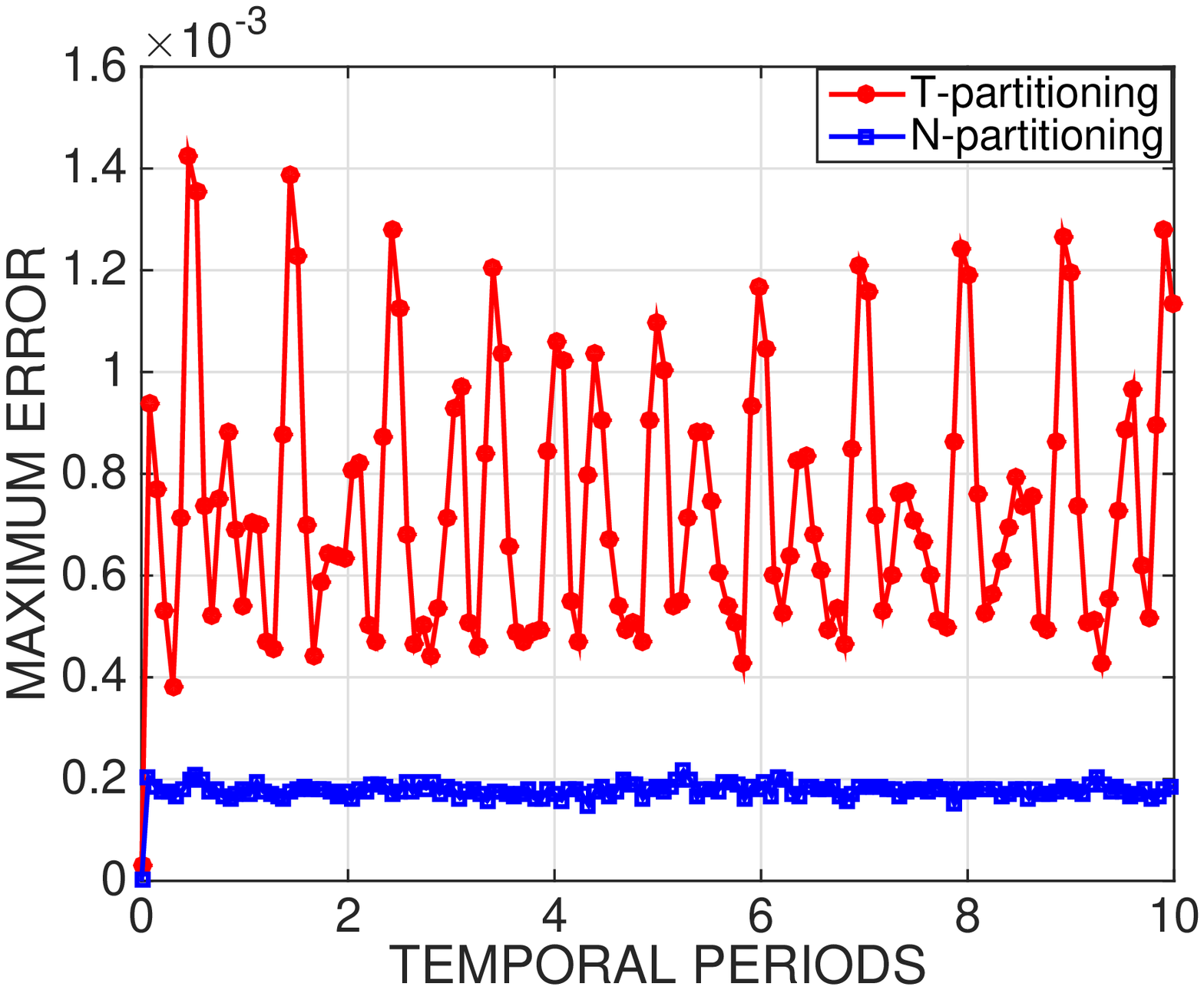}\label{fig:SixthCavity}}
	\caption{Maximum error of the numerical experiment with the circular cavity}
\end{figure}
We employ the fourth and sixth order SBP--SAT method to propagate the wave for ten periods, and show the recorded maximum errors in Figure \ref{fig:FourthCavity} and \ref{fig:SixthCavity}. In both cases, the maximum error with the T--partitioning is about three times larger than that with the N--partitioning. The is not surprising because the mesh with the T--partitioning has less grid points than the mesh with the N--partitioning, and the corresponding scheme with the T--partitioning has one order lower convergence rate than that of the N--partitioning. It does not seem to improve the efficiency by using T--junction interfaces for this case.

Although using T--junction interfaces introduces a larger error in the solution, it could be beneficial for a problem with a more complex geometry. For example, if there are several cavities in the domain, an N--partitioning that only allows conforming mesh blocks would produce a large number of small mesh blocks. With a higher order accurate scheme, the boundary stencil gets wider and the number of grid points in each direction must be large enough in every mesh block. It is therefore over--resolved in those small mesh blocks and results in a suboptimal performance of the numerical scheme, and T--junction interfaces could be desirable.

\subsubsection{A circular inclusion}
Consider a circular domain of radius $a=1$ embedded in an infinite surrounding medium of differing material with wave speed $c$. Let the wave speed $c'$ of the circular domain be such that $c'<c$ and let an incoming time--harmonic plane wave $u_I(x,y,t) = e^{i(\omega t - \gamma x)}, \gamma=\frac{\omega}{c}$ travel in the $x$--direction and impinge on the circular inclusion. The resulting field consists of the incoming wave $u_I$, as well as the scattered and diffracted waves $u_S$ and $u_D$, respectively. The conditions at the interface of the circular inclusion are
\begin{equation}
	\label{eq:myPart4}
	\begin{split}
		&u_I + u_S = u_D,\\
		&c \frac{\partial}{\partial n}\left(u_I + u_S\right) = c' \frac{\partial}{\partial n} u_D, 
	\end{split} \text{on}\ x^2 + y^2 = 1,
\end{equation}
where $\frac{\partial}{\partial n}$ denotes the normal derivative on the interface. Since $c'<c$, the short wavelength occurs inside the circular domain. An analytical expression for the solution is given in \cite[pp.~667]{Balanis1989}.

\begin{figure}
        \includegraphics[trim=0cm 3.5cm 0cm 2.5cm, clip=true, width= \textwidth]{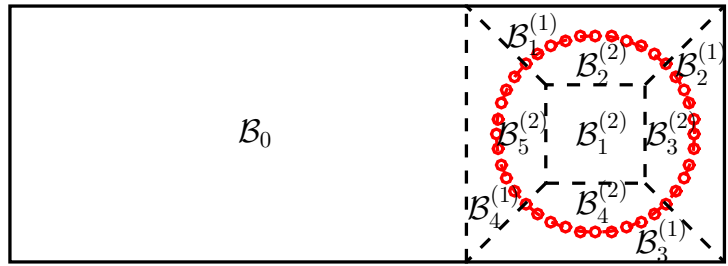} 
        	\caption{Computational domain of the circular inclusion experiment}
	\label{fig:Domain_Inclusion}
\end{figure}

In the numerical experiments we take $\omega = 2 \pi$, $c=1$ and $c' = 1/10$, which give a wavelength of $1$ and $1/10$ outside and inside the circular inclusion, respectively. To resolve the geometric features, the computational domain is decomposed into 10 conforming blocks as shown in Figure \ref{fig:Domain_Inclusion}. We take the side length $2D = 2.6$ for the square block outside the circular inclusion, and the side length $2d = 0.7\sqrt{2}$ for the square block inside the circular inclusion. Both square blocks are centered at the origin. The Cartesian block $[-5.9,-1.3]\times[-1.3,1.3]$ is then attached to the left of this representation. Firstly, we only use conforming grid interfaces with the numbers of grid points in each block given in Table \ref{CBlock_table}. The resolution outside the circular inclusion is about $16$ and $51$ points per wavelength in the horizontal and vertical direction, respectively.  Inside the circular inclusion the waves are resolved by about $10$ grid points per wavelength in both directions. Hence, the waves are significantly over--resolved in the vertical direction outside the inclusion, which leads to a suboptimal efficiency of the numerical scheme. 

To amend the over--resolution, we partition the computational domain in the same way as above but use non--conforming interfaces denoted by the small red circles in Figure \ref{fig:Domain_Inclusion}. The non--conforming grid interfaces are handled by the interpolation and projection operators. The numbers of grid points in each block are chosen as in Table \ref{NCBlock_table}.  Now the resolution is reduced to about $26$ grid points per wavelength in the vertical direction outside the circular inclusion. The interface conditions \eqref{eq:myPart4} are imposed numerically with the SAT technique and at outer boundaries the exact solution is injected at all times. In \cite{Duru2014}, the SBP finite difference method applied to the wave equation with the injection method to impose the Dirichlet boundary condition is proved to be stable. 

\begin{table}
\begin{minipage}{0.48\textwidth}
\centering
\begin{tabular}{c cc}
\toprule
Block & $N_\xi$ & $N_\eta$ \\ \midrule
 $\mathcal{B}_0$  & $51$ & $101$  \\ \midrule
 $\mathcal{B}^{(1)}_1$--$ \mathcal{B}^{(1)}_4$  & $101$    &  $26$ \\ \midrule
 $\mathcal{B}^{(2)}_1$--$ \mathcal{B}^{(2)}_4$  &  $101$   &  $51$ \\  \midrule
 $\mathcal{B}^{(2)}_5$  & $101$    &  $101$ \\ \bottomrule
\end{tabular}
\caption{Number of grid points with conforming grid interfaces}
\label{CBlock_table}
\end{minipage}
\hfill
\begin{minipage}{0.48\textwidth}
\centering
\begin{tabular}{c cc}
\toprule
Block & $N_\xi$ & $N_\eta$ \\ \midrule
 $\mathcal{B}_0$  & $51$ & $51$  \\ \midrule
 $\mathcal{B}^{(1)}_1$ --$ \mathcal{B}^{(1)}_4$  & $51$    &  $26$ \\ \midrule
 $\mathcal{B}^{(2)}_1$ -- $\mathcal{B}^{(2)}_4$  &  $101$   &  $51$ \\  \midrule
 $\mathcal{B}^{(2)}_5$  & $101$    &  $101$ \\ \bottomrule
\end{tabular}
\caption{Number of grid points with non--conforming grid interfaces}
\label{NCBlock_table}
\end{minipage}
\end{table}

The solution is propagated numerically for 10 temporal periods by the SBP--SAT method and the relative maximum error is recorded at each time step. 
In Figures \ref{fig:4th_inclusion} and \ref{fig:6th_inclusion} we plot the recorded relative maximum error as functions of time. Here we see that the errors are similar in both cases. The grid with conforming interfaces has $46460$ grid points, whereas the grid with non--conforming interfaces has $38710$ grid points. The smallest grid size is determined by the resolution inside the circular inclusion, for this reason the time step $\Delta t = 4 \times 10^{-4}$ for the sixth order SBP-SAT method and $\Delta t = 5 \times 10^{-4}$ for the fourth order SBP-SAT method are the same for both grids. We conclude that even though the formal order of accuracy is lowered by using blocks with non-conforming interfaces it can be a beneficial strategy within the SBP-SAT framework when the physical domain contains regions that require a higher density of grid points. We also note that for more complex multi--block domains consisting of a larger number of blocks the benefits of using blocks with non--conforming grid interfaces are expected to increase. 

In the experiment with the sixth order SBP--SAT scheme with the interpolation operators, the numerical solution blows up quickly, which indicates that the scheme is unstable. The corresponding scheme with the projection operators is stable. According to the stability analysis in \S\ref{sec-stability}, for the sixth order accurate scheme condition (\ref{inequality}) holds for neither the interpolation operator nor the projection operator. If the smallest eigenvalue of $\Xi_{L/R}$ is non--negative, then an energy estimate exists that ensures stability. The smallest eigenvalue of $\Xi_{L/R}$ scaled by the mesh size is in the magnitude of $-10^{-1}$ for the sixth order accurate interpolation operators, and $-10^{-5}$ for the sixth order accurate projection operators. The violation of (\ref{inequality}) is much weaker for the projection operator than for the interpolation operator, which explains the observation in the numerical experiments. 

\begin{figure}
	\subfloat[Fourth order accurate scheme]{\includegraphics[width=0.49\textwidth]{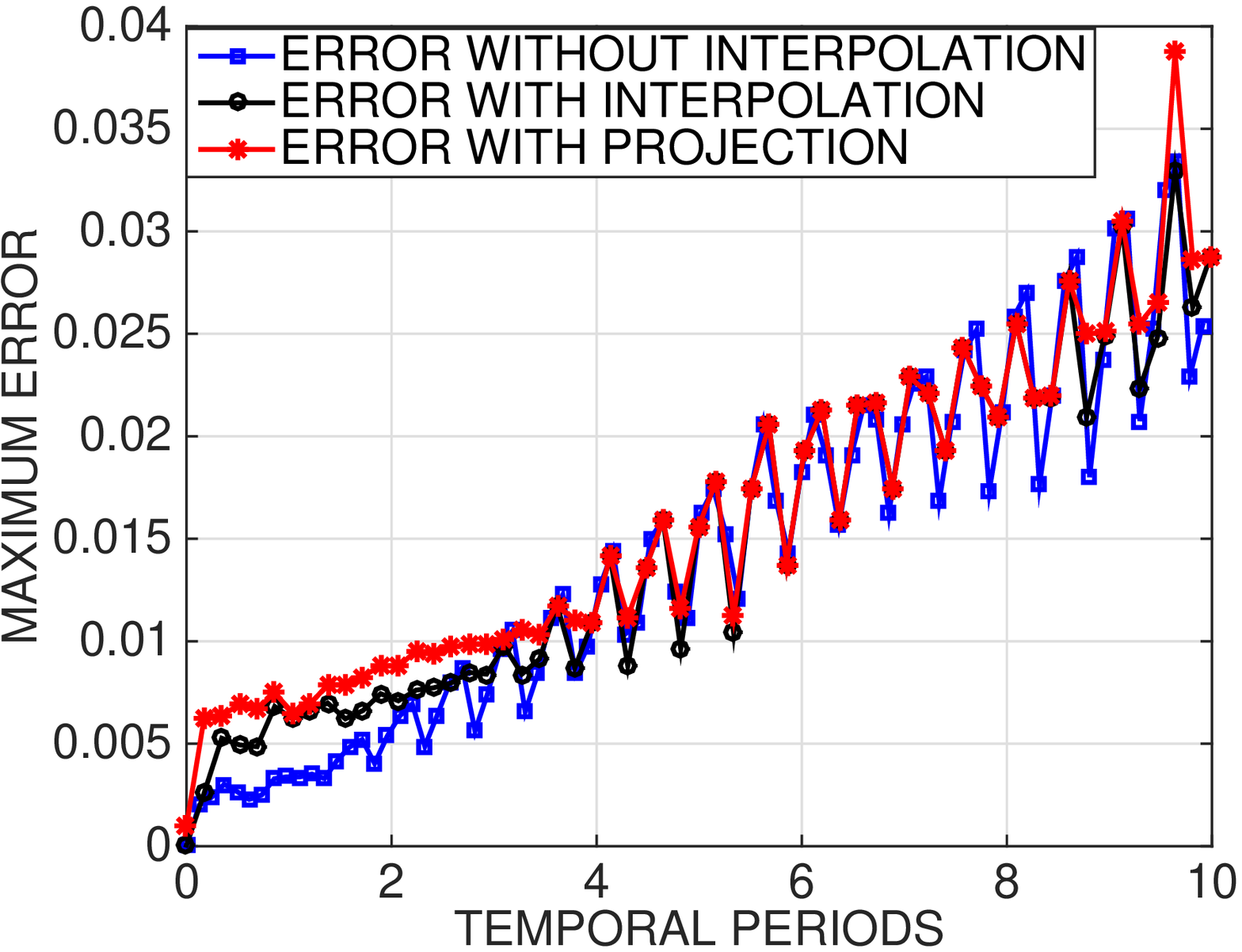}\label{fig:4th_inclusion}}
	\subfloat[Sixth order accurate scheme]{\includegraphics[width=0.49\textwidth]{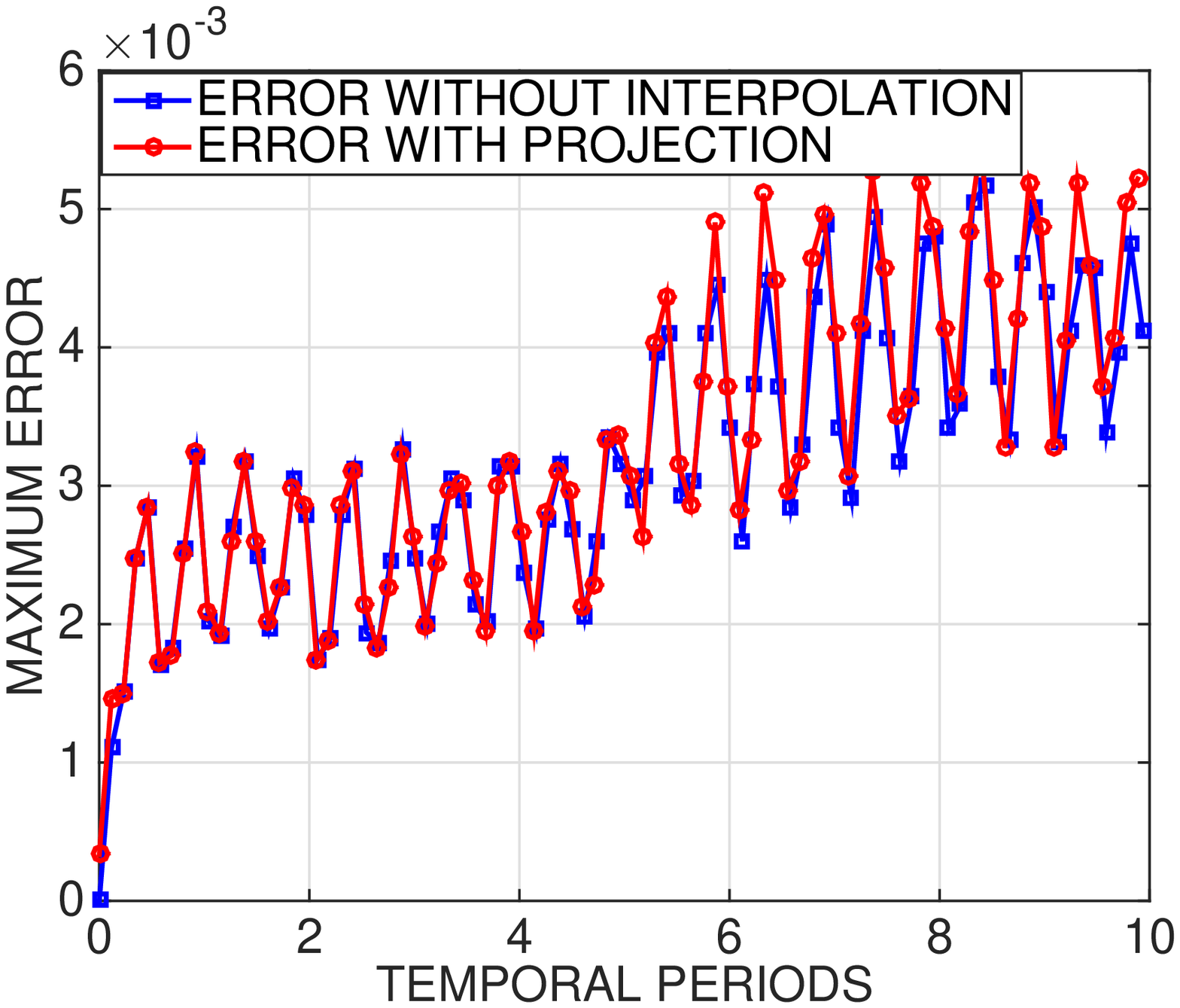}\label{fig:6th_inclusion}}
	\caption{Maximum error of the numerical experiment with the circular inclusion}
\end{figure}

\section{Conclusion and outlook}\label{sec-Conclusion}
In this work, high order accurate SBP finite difference operators are used to discretize the wave equation in the second order form on a block--structured mesh. Adjacent mesh blocks are patched together by imposing suitable interface conditions via the SAT technique. The emphasis is placed on the numerical treatment of non--conforming grid interfaces by the interpolation operators and projection operators, which are also compared in terms of the stability and accuracy properties. In contrast to first order hyperbolic equation, stability of the scheme for the second order wave equation introduces an extra condition on the numerical treatment of non--conforming grid interfaces. This condition is satisfied for the second and fourth order accurate cases, and an energy estimate is derived to ensure stability. For higher order accurate schemes, the extra stability condition is violated. We show by the eigenvalue analysis that the violation is stronger with interpolation operators than with projection operators. Unphysical growths are observed in the numerical experiments with high order interpolation operators, whereas with projection operators the scheme is stable in all the experiments we have conducted. 

We have also performed a truncation error analysis and an investigation of the convergence property for the scheme, which indicates that the convergence rate is one order lower than that in the corresponding case with conforming grid interfaces. This is verified in numerical experiments. The efficiency studies show that for a problem with a very complex geometry, it could be beneficial to use non--conforming grid interfaces. 

For high order accurate interpolation operators and projection operators, there are free parameters left in the construction process. It is desirable to tune the free parameters so that the extra stability condition is satisfied. However, the resulting nonlinear problem seems difficult to solve. To overcome the accuracy reduction, more accurate interpolation operators or new ways of imposing interface conditions are needed.  

\section*{Appendix}
We give a detailed description for the construction of the two grids used in the efficiency studies in \S\ref{sec:efficiency_study}. 
A general block $\mathcal{B}_i$ of a decomposition has four boundaries defined by the parametrized curves 
\begin{equation*}
	\begin{split}
		&\mathcal{C}_{iS} = \left(x_{iS}(\xi),y_{iS}(\xi)\right),\quad \mathcal{C}_{iN} = \left(x_{iN}(\xi),y_{iN}(\xi)\right),\\
		&\mathcal{C}_{iW} = \left(x_{iW}(\eta),y_{iW}(\eta)\right), \quad \mathcal{C}_{iE} = \left(x_{iE}(\eta),y_{iE}(\eta)\right),
	\end{split} 
	\end{equation*}   
where $0\leq \xi \leq 1, 0\leq \eta \leq 1$. $\mathcal{C}_{iS}$ and $\mathcal{C}_{iN}$ describe one pair of opposing sides and $\mathcal{C}_{iW}$ and $\mathcal{C}_{iE}$ the other pair. Let $\mathcal{P}_{iSW}$ denote the point of intersection between the curves $\mathcal{C}_{iS}$ and $\mathcal{C}_{iW}$ e.t.c. A bijection $(x,y) = \mathcal{T}_i(\xi,\eta)$ from the unit square $\mathcal{S} = [0,1]^2$ to the block $\mathcal{B}_i$ of the decomposition is given by the transfinite interpolation \cite{Knupp_Fundamentals_of_grid_generation} as
\begin{equation*}
	\begin{split}
		\mathcal{T}_i(\xi,\eta) &= (1-\eta) \mathcal{C}_{iS} + \eta \mathcal{C}_{iN} + (1-\xi) \mathcal{C}_{iW} + \xi \mathcal{C}_{iE}\\
									&-\xi \eta \mathcal{P}_{iNE} - \xi (1-\eta) \mathcal{P}_{iSE} - \eta (1-\xi) \mathcal{P}_{iNW} - (1-\xi)(1-\eta)\mathcal{P}_{iSW}.
	\end{split}
\end{equation*} 
The unit square $\mathcal{S}$ is discretized by the points
\begin{equation*}
	\label{eq:discretized_unit_square}
	\begin{split}
		&\xi_j = j h_{\xi}, h_{\xi} = 1/(N_{i\xi}-1), j = 0, \dots, N_{i\xi}-1,\\
		&\eta_k = k h_{\eta}, h_{\eta} = 1/(N_{i\eta}-1), k = 0, \dots, N_{i\eta}-1,
	\end{split}
\end{equation*}
where $N_{i\xi}$ and $N_{i\eta}$ are integers determining the number of grid points in the spatial directions of the discretization of the block $\mathcal{B}_i$. The corresponding grid points are computed as 
\begin{equation*}
	\label{eq:compute_grid_points}
	(x_j,y_k) = \mathcal{T}_i(\xi_j,\eta_k).
\end{equation*}
We now give details on how the grids in $\mathcal{D}_1 - \mathcal{D}_2$ are constructed. $\mathcal{D}_1$ represents a circular cavity in an infinite surrounding media. We construct the computational grid such that the cavity of radius $a$ is centered inside a square of side $2D$. This is done by introducing four blocks $\mathcal{B}^{(1)}_1$--$\mathcal{B}^{(1)}_4$. The bounding curves of block $\mathcal{B}^{(1)}_1$ are given by
\begin{equation*}
	\label{eq:D1}
	\begin{split}
		&\mathcal{C}^{(1)}_{1S} = a\left(\xi \sqrt{2} - 1/\sqrt{2}, \sqrt{1-(\xi \sqrt{2} - 1/\sqrt{2})^2}\right),\\
		&\mathcal{C}^{(1)}_{1N} = \left(D, 2D \xi - D\right),\\
		&\mathcal{C}^{(1)}_{1W} = \left(-\eta (D-a/\sqrt{2}) - a/\sqrt{2}, \eta (D-a/\sqrt{2}) + a/\sqrt{2}\right),\\
		&\mathcal{C}^{(1)}_{1E} = \left( \eta (D-a/\sqrt{2}) + a/\sqrt{2}, \eta (D-a/\sqrt{2}) + a/\sqrt{2}\right).
	\end{split} 
\end{equation*} 
The bounding curves of the remaining blocks $\mathcal{B}^{(1)}_2$--$\mathcal{B}^{(1)}_4$ that constitute the square with the cavity at the center are obtained via rotation by a factor $\pi/2$,
\begin{equation}
	\label{eq:rot_pi_div_2}
	\mathcal{C}^{(1)}_{ij} = \mathcal{C}^{(1)}_{i-1j} \begin{bmatrix} \cos{\pi/2} & -\sin{\pi/2}\\ \sin{\pi/2} & \cos{\pi/2}\end{bmatrix},\ i = 2,3,4,\ j = S,N,W,E.
\end{equation} 
The domain $\mathcal{D}_1$ can now be represented by attaching the grid representing the cavity to one or more Cartesian blocks.

The circular inclusion of $\mathcal{D}_2$ is decomposed into five blocks $\mathcal{B}^{(2)}_1$--$\mathcal{B}^{(2)}_5$. The block $\mathcal{B}^{(2)}_1$ is a square at the center of the circular inclusion with corners at the points $(\pm ad,\pm ad), 0 < d < \sqrt{2}/2$ defined by the bounding curves,
\begin{equation*}
	\label{eq:D21}
	\begin{split}
		&\mathcal{C}^{(2)}_{1S} = a\left(2d \xi-d, -d\right),\ \mathcal{C}^{(2)}_{1N} = a\left(2d\xi-d, d\right),\\
		&\mathcal{C}^{(2)}_{1W} = a\left(-d, 2d\eta-d\right),\ \mathcal{C}^{(2)}_{1E} = a\left(d, 2d\eta-d\right).
	\end{split} 
\end{equation*} 
Here $a$ is the radius of the circular inclusion. The block $\mathcal{B}^{(2)}_2$ is defined by its bounding curves, 
\begin{equation*}
	\label{eq:D22}
	\begin{split}
		&\mathcal{C}^{(2)}_{2S} = \mathcal{C}^{(2)}_{1N},\\
		&\mathcal{C}^{(2)}_{2N} = \mathcal{C}^{(1)}_{1S},\\
		&\mathcal{C}^{(2)}_{2W} = a\left(-\eta (\sqrt{2}/2-d)-d, \eta(\sqrt{2}/2-d) + d\right),\\
		&\mathcal{C}^{(2)}_{2E} = a\left( \eta (\sqrt{2}/2-d) + d, \eta (\sqrt{2}/2-d) + d\right).
	\end{split} 
\end{equation*} 
The bounding curves of the remaining blocks $\mathcal{B}^{(2)}_3$--$\mathcal{B}^{(2)}_5$ that constitute the circular inclusion of $\mathcal{D}_2$ are obtained via rotations by a factor $\pi/2$ as in \eqref{eq:rot_pi_div_2}. The inclusion is then centered inside a square of side $2D$ by attaching it to the blocks $\mathcal{B}^{(1)}_1$--$\mathcal{B}^{(1)}_4$ above. The circular inclusion is now the union of the nine blocks $\mathcal{B}^{(1)}_1$--$\mathcal{B}^{(1)}_4$ and $\mathcal{B}^{(2)}_1$--$\mathcal{B}^{(2)}_5$. The domain $\mathcal{D}_2$ can now be represented by attaching the grid representing the inclusion to one or more Cartesian blocks.


\begin{thebibliography}{1}

\bibitem{Balanis1989}{\sc C.~A. Balanis},
{\em Advanced Engineering Electromagnetics},
Wiley, 1989.

\bibitem{Carpenter1994} {\sc M.~H. Carpenter, D. Gottlieb and S. Abarbanel},
{\em Time-stable boundary conditions for finite-difference schemes solving hyperbolic systems: methodology and application to high-order compact schemes},
J. Comput. Phys., 111(1994), pp.~220--236.

\bibitem{Fernandez2014}{\sc D.~C. Del Rey Fern\'{a}ndez, J.~E. Hicken and D.~W. Zingg},
{\em Review of summation-by-parts operators with simultaneous approximation terms for the numerical solution of partial differential equations},
Comput. \& Fluids, 95(2014), pp.~171--196.

\bibitem{Duru2014}{\sc K. Duru, G. Kreiss and K. Mattsson},
{\em Stable and high--order accurate boundary treatments for the elastic wave equation on second--order form},
SIAM J. Sci. Comput., 36(2014), pp.~A2787--A2818.

\bibitem{Gassner2013}{\sc G.~J. Gassner},
{\em A skew--symmetric discontinuous Galerkin spectral element discretization and its relation to SBP--SAT finite difference methods},
SIAM J. Sci. Comput., 35(2013), pp.~A1233--A1253.

\bibitem{Graff1991}{\sc K.~F. Graff},
{Wave Motion in Elastic Solids}, Dover Publications, 1991.

\bibitem{Gustafsson2008}{\sc B. Gustafsson},
{\em High Order Difference Methods for Time Dependent PDE}, Springer--Verlag, Berlin Heidelberg, 2008.

\bibitem{Gustafsson2013}{\sc B. Gustafsson, H.~O. Kreiss and J. Oliger},
{\em Time-Dependent Problems and Difference Methods}, Wiley, New Jersey, 2013.

\bibitem{Knupp_Fundamentals_of_grid_generation} {\sc P. Knupp and S. Steinberg}, 
{\em Fundamentals of Grid Generation}, CRC Press, 1993. 
 
\bibitem{Kozdon2015}{\sc J.~E. Kozdon and L.~C. Wilcox},
{\em Provably stable, general purpose projection operators for high--order finite difference methods}, arXiv:1410.5746v2 [math.NA].
 
\bibitem{Kreiss1972}{\sc H.~O. Kreiss and J. Oliger},
{\em Comparison of accurate methods for the integration of hyperbolic equations}, Tellus XXIV, 24(1972), pp.~199--215.

\bibitem{Kreiss2002}{\sc H.~O. Kreiss, N.~A. Petersson and J. Ystr\"{o}m},
{\em Difference approximations for the second order wave equation},
SIAM J. Numer. Anal., 40(2002), pp.~1940--1967.

\bibitem{Kreiss1974}{\sc H.~O. Kreiss and G. Scherer},
{\em Finite element and finite difference methods for hyperbolic partial differential equations},
Mathematical aspects of finite elements in partial differential equations, Symposium proceedings (1974), pp.~195--212.

\bibitem{Mattsson2012}{\sc K. Mattsson},
{\em Summation by parts operators for finite difference approximations of second--derivatives with variable coefficients}, J. Sci. Comput., 51(2012), pp.~650--682.

\bibitem{Mattsson2013} {\sc K. Mattsson and M. Almquist},
{\em A solution to the stability issues with block norm summation by parts operators}, J. Comput. Phys., 253(2013), pp.~418--442.

\bibitem{Mattsson2010} {\sc K. Mattsson and M.~H. Carpenter},
{\em Stable and accurate interpolation operators for high-order multiblock finite difference methods}, SIAM J. Sci. Comput., 32(2010), pp.~2298--2320.

\bibitem{Mattsson2008} {\sc K. Mattsson, F. Ham and G. Iaccarino}, 
{\em Stable and accurate wave-propagation in discontinuous media}, J. Comput. Phys., 227(2008), pp.~8753--8767.

\bibitem{Mattsson2009} {\sc K. Mattsson, F. Ham and G. Iaccarino},
{\em Stable boundary treatment for the wave equation on second-order form}, J. Sci. Comput., 41(2009), pp.~366--383.

\bibitem{Mattsson2004} {\sc K. Mattsson and J. Nordstr\"{o}m},
{\em Summation by parts operators for finite difference approximations of second derivatives}, J. Comput. Phys., 199(2004), pp.~503--540.

\bibitem{Mattsson2008b}{\sc K. Mattsson, M. Sv\"{a}rd and M. Shoeybi},
{\em Stable and accurate schemes for the compressible Navier--Stokes equations},
J. Comput. Phys., 227(2008), pp.~2293--2316.

\bibitem{Nissen2015}{\sc A. Nissen, K. Kormann, M. Grandin and K. Virta},
{\em Stable difference methods for block-oriented adaptive grids}, J. Sci. Comput, DOI: 10.1007/s10915-014-9969-z.

\bibitem{Nissen2012} {\sc A. Nissen, G. Kreiss and M. Gerritsen},
{\em Stability at non--conforming grid interfaces for a high order discretization of the Schr\"{o}dinger equation}, J. Sci. Comput., 53(2012), pp.~528--551.

\bibitem{Nissen2013} {\sc A. Nissen, G. Kreiss and M. Gerritsen},
{\em High order stable finite difference methods for the Schr\"{o}dinger equation}, J. Sci. Comput., 55(2013), pp.~173--199.

\bibitem{Petersson2010}{\sc N.~A. Petersson and B. Sj\"{o}green},
{\em Stable grid refinement and singular source discretization for seismic wave simulations}, Commun. Comput. Phys., 8(2010), pp.~1074--1110. 

\bibitem{Strand1994}{\sc B. Strand},
{\em Summation by parts for finite difference approximations for d/dx}, J. Comput. Phys. 110(1994), pp.~47--67.

\bibitem{Svard2014}{\sc M. Sv\"{a}rd and J. Nordstr\"{o}m},
{\em Review of summation-by-parts schemes for initial-boundary-value problems}, J. Comput. Phys., 268(2014), pp.~17--38.

\bibitem{Virta2014}{\sc K. Virta and K. Mattsson},
{\em Acoustic wave propagation in complicated geometries and heterogeneous media}, J. Sci. Comput., 61(2014), pp.~90--118.

\bibitem{Wang2015}{\sc S. Wang and G. Kreiss},
{\em Convergence of summation--by--parts finite difference methods for the wave equation}, arXiv:1503.08926v2 [math.NA].

\end{thebibliography}
\end{document}